# SERIAL AND NONSERIAL SIGN-AND-RANK STATISTICS: ASYMPTOTIC REPRESENTATION AND ASYMPTOTIC NORMALITY

By Marc Hallin,[1] Catherine Vermandele and Bas Werker

*Université Libre de Bruxelles, Université Libre de Bruxelles and Tilburg University*

The classical theory of rank-based inference is entirely based either on *ordinary* ranks, which do not allow for considering location (intercept) parameters, or on signed ranks, which require an assumption of symmetry. If the median, in the absence of a symmetry assumption, is considered as a location parameter, the maximal invariance property of ordinary ranks is lost to the ranks and the signs. This new maximal invariant thus suggests a new class of statistics, based on ordinary ranks *and* signs. An asymptotic representation theory à la Hájek is developed here for such statistics, both in the nonserial and in the serial case. The corresponding asymptotic normality results clearly show how the signs add a separate contribution to the asymptotic variance, hence, potentially, to asymptotic efficiency. As shown by Hallin and Werker [*Bernoulli* **9** (2003) 137–165], conditioning in an appropriate way on the maximal invariant potentially even leads to semiparametrically efficient inference. Applications to semiparametric inference in regression and time series models with median restrictions are treated in detail in an upcoming companion paper.

**1. Introduction.** The classical theory of rank-based inference is entirely based either on *ordinary ranks* or on *signed ranks*. Ranks indeed are maximal invariant with respect to the group of continuous order-preserving transformations, a group that generates the null hypothesis of absolutely continuous independent white noise (no location restriction), whereas signed ranks (i.e., the signs along with the ranks of absolute values) are maximal invariant

Received February 2003; revised August 2004.
[1]Supported by an I.A.P. grant of the Belgian Federal Government and an Action de Recherche Concertée from the Communauté française de Belgique.

under the subgroup that generates the subhypothesis of symmetric (with respect to the origin) independent white noise.

Now, in most statistical models a location parameter for the error term is usually specified to be zero: regression and analysis of variance models, stationary autoregressive moving average (ARMA) models and so on. Symmetric white noise allows for such an identification, at the expense, however, of a symmetry assumption that in practice is often quite unrealistic. In addition, the trouble with independent white noise without further restrictions is that it does not allow for identifying any location parameter.

This location parameter in most applied work is the mean—a heritage of Gaussian models—but could be the median as well. Zero-median noise is certainly as natural as zero-mean noise. In a semiparametric context, it is even more satisfactory, because it does not require any moment assumption on the densities under consideration. Median regression and autoregression models have, therefore, recently attracted much attention: see, for instance, [12, 14, 15, 17, 21], to quote only a few. Moreover, from the point of view of statistical inference, the assumption of zero-median noise is also more convenient, since it induces more structure. The hypothesis of zero-mean white noise indeed is not invariant under any nontrivial group of transformations, so group invariance arguments cannot be invoked in models that involve zero-mean noise. The situation is quite different for the hypothesis of zero-median noise, which is generated by the group of all continuous order-preserving transformations $g$ such that $g(0) = 0$. A maximal invariant for this group is the vector of ordinary ranks, along with the vector of signs. Hallin and Werker [11] have shown that, in such a situation, semiparametric efficiency is achieved by conditioning with respect to a *maximal* invariant. Maximality of the invariant here is essential: conditioning, for example, on the ranks when the signs and ranks, not the ranks alone, are maximal invariant generally induces an avoidable loss of efficiency.

Invariance and semiparametric efficiency arguments in such models thus lead to the new concept of sign-and-rank-based statistics, which involve both signs and ranks. This new concept is more natural than the traditional rank-based one in all models that include a location (intercept) parameter, but also in models such as stationary ARMA models, where the noise is inherently centered. The objective of the present paper is a detailed study of the class of linear sign-and-rank statistics for which we provide Hájek-type asymptotic representation and asymptotic normality results. These results readily allow for building new rank-based tests for a variety of problems in one-, two- and $k$-sample location, regression, ARMA and related models without making any symmetry assumptions on the underlying error densities. They also form a basis for the construction of semiparametrically efficient procedures in median constrained models (see [10]).

The paper is organized as follows. Section 2 briefly introduces several concepts of white noise: *independent*, *independent with zero mean*, *independent with zero median* and *independent symmetric* white noises. We recall how the invariance principle for each of these concepts, but for white noise with zero mean, leads to a different concept of ranks and/or signs—the right concept for median-centered white noise being the signs and ranks. Sections 3 and 4 propose a systematic investigation of (linear) nonserial and serial sign-and-rank statistics. These new statistics, which are measurable with respect to the vectors of ranks and signs, are studied along the same lines as the classical linear rank statistics (see, e.g., [3] for the nonserial context; see [5] and [7] for the serial context) and the linear signed-rank statistics (see [3] and [13] for the nonserial context; see [7] for the serial context). However, the nonindependence between the ranks and the signs (in sharp contrast with the traditional context of signed ranks, where the signs and the ranks of absolute values are mutually independent) requires a more delicate treatment. Section 5 concludes with an empirical study: simulations very clearly show that the proposed procedures quite significantly outperform their classical counterparts based on either parametric correlograms or traditional ranks—the more skewed the underlying densities, the more significant the efficiency gain.

## 2. White noise and group invariance.

2.1. *White noise and semiparametric statistical models.* Whatever the concept of ranks, rank-based inference applies in the context of semiparametric models under which the distribution of some observed $n$-tuple $\mathbf{Y}^{(n)} := (Y_1^{(n)}, \ldots, Y_n^{(n)})'$ belongs to a family of distributions of the form

$$\{\mathrm{P}^{(n)}_{f;\boldsymbol{\theta}}, \boldsymbol{\theta} \in \boldsymbol{\Theta} \subseteq \mathbb{R}^K, f \in \mathcal{F}\}, \tag{2.1}$$

where $\boldsymbol{\theta}$ denotes some finite-dimensional parameter of interest and $f$ denotes some unspecified density (densities throughout are tacitly taken with respect to the Lebesgue measure over the real line) that plays the role of a nonparametric nuisance. This distribution $\mathrm{P}^{(n)}_{f;\boldsymbol{\theta}}$, in general, is described by means of (i) a *residual function*, namely, a family $\{\mathfrak{Z}^{(n)}_{\boldsymbol{\theta}}, \boldsymbol{\theta} \in \boldsymbol{\Theta}\}$ of invertible functions indexed by $n$ and $\boldsymbol{\theta}$ that map the observation $\mathbf{Y}^{(n)}$ onto an $n$-tuple of *residuals*

$$\mathfrak{Z}^{(n)}_{\boldsymbol{\theta}}(\mathbf{Y}^{(n)}) = \mathbf{Z}^{(n)}(\boldsymbol{\theta}) := (Z_1^{(n)}(\boldsymbol{\theta}), \ldots, Z_n^{(n)}(\boldsymbol{\theta}))',$$

and (ii) a concept of *white noise* with (marginal) density $f$ such that $\mathbf{Y}^{(n)}$ has distribution $\mathrm{P}^{(n)}_{f;\boldsymbol{\theta}}$ iff $\mathbf{Z}^{(n)}(\boldsymbol{\theta})$ is white noise with (marginal) density $f$.

We concentrate on four particular forms of white noise. Define $\mathcal{F} := \{f : f(x) > 0, x \in \mathbb{R}\}$ as the set of all nonvanishing densities over the real line, let $\mathcal{F}_* := \{f \in \mathcal{F} : \int_{-\infty}^{\infty} z f(z)\,dz = 0\}$ be the subset of all densities in $\mathcal{F}$ with zero mean, let $\mathcal{F}_0 := \{f \in \mathcal{F} : \int_{-\infty}^{0} f(z)\,dz = \int_0^\infty f(z)\,dz = 1/2\}$ be the set of densities in $\mathcal{F}$ having zero median and let $\mathcal{F}_+ := \{f \in \mathcal{F} : f(-z) = f(z), z \in \mathbb{R}\}$ be the set of densities in $\mathcal{F}$ that are symmetric with respect to the origin. Denote the following terms:

(a) *Independent white noise*: Let $\mathcal{H}_f^{(n)}$ denote the hypothesis under which the random vector $\mathbf{Z}^{(n)} = (Z_1^{(n)}, \ldots, Z_n^{(n)})'$ is a realization of length $n$ of an *independent* white noise; that is, $Z_i^{(n)}$, $i = 1, \ldots, n$, are i.i.d. with density $f \in \mathcal{F}$.

(b) *Zero-mean independent white noise*: Let $\mathcal{H}_{*;f}^{(n)}$ denote the hypothesis under which $\mathbf{Z}^{(n)}$ is a realization of length $n$ of an *independent with zero-mean* white noise; that is, $Z_i^{(n)}$, $i = 1, \ldots, n$, are i.i.d. with density $f \in \mathcal{F}_*$.

(c) *Zero-median independent white noise*: Let $\mathcal{H}_{0;f}^{(n)}$ denote the hypothesis under which $\mathbf{Z}^{(n)}$ is a realization of length $n$ of an *independent with zero-median* white noise; that is, $Z_i^{(n)}$, $i = 1, \ldots, n$, are i.i.d. with density $f \in \mathcal{F}_0$.

(d) *Symmetric independent white noise*: Let $\mathcal{H}_{+;f}^{(n)}$ denote the hypothesis under which $\mathbf{Z}^{(n)}$ is a realization of length $n$ of an *independent symmetric* white noise; that is, $Z_i^{(n)}$, $i = 1, \ldots, n$, are i.i.d. with density $f \in \mathcal{F}_+$.

The notation $\mathcal{H}^{(n)}$, $\mathcal{H}_*^{(n)}$, $\mathcal{H}_0^{(n)}$ and $\mathcal{H}_+^{(n)}$ is used whenever the underlying density function $f$ remains unspecified within $\mathcal{F}$, $\mathcal{F}_*$, $\mathcal{F}_0$ or $\mathcal{F}_+$, respectively. In practice, of course, the role of the random variables $Z_i^{(n)}$ is played by the residuals $Z_i^{(n)}(\boldsymbol{\theta})$ $(i = 1, \ldots, n)$ associated with a specific value $\boldsymbol{\theta}$ of the parameter in the statistical model under consideration.

The independent white noise hypothesis $\mathcal{H}^{(n)}$ is most general, but does not allow for identifying location parameters. A classical attitude, when location is to be identified, consists in assuming that the underlying white noise density has zero mean, that is, adopting $\mathcal{H}_*^{(n)}$. As already explained, an often-used alternative solution requires the median (instead of the mean) of the white noise density to be zero, leading to $\mathcal{H}_0^{(n)}$. The additional assumption of symmetry yields $\mathcal{H}_+^{(n)}$.

2.2. *Group invariance*: *ranks*, *signed ranks*, *and signs and ranks.* Let $\mathcal{E}^{(n)} := (\mathbb{R}^n, \mathcal{B}^n, \mathcal{P}^{(n)} := \{\mathrm{P}_{\boldsymbol{\theta};f}^{(n)}, \boldsymbol{\theta} \in \boldsymbol{\Theta}, f \in \mathcal{F}\})$ be characterized (in the sense of Section 2.1) by the residual function $\mathfrak{Z}_{\boldsymbol{\theta}}^{(n)}$ and the white noise concept $\mathcal{H}^{(n)}$. Denote by $\mathfrak{G}$ the set of all continuous, strictly monotone increasing functions $g : \mathbb{R} \to \mathbb{R}$ such that $\lim_{x \to \pm\infty} g(x) = \pm\infty$, define $\mathcal{G}_g^{(n)} : \mathbf{z} = (z_1, \ldots, z_n)' \in$ $\mathbb{R}^n \mapsto \mathcal{G}_g^{(n)}(\mathbf{z}) := (g(z_1), \dots, g(z_n))' \in \mathbb{R}^n$ and consider the group (acting on $\mathbb{R}^n$)

$$\mathcal{G}_{\boldsymbol{\theta}}^{(n)}, \circ := \{(\mathfrak{Z}_{\boldsymbol{\theta}}^{(n)})^{-1} \circ \mathcal{G}_g^{(n)} \circ \mathfrak{Z}_{\boldsymbol{\theta}}^{(n)}, g \in \mathfrak{G}\}, \circ.$$

This group (known as the group of *order-preserving transformations* of residuals) clearly is a generating group for the fixed-$\boldsymbol{\theta}$ submodel $\mathcal{E}^{(n)}(\boldsymbol{\theta}) := (\mathbb{R}^n, \mathcal{B}^n, \mathcal{P}^{(n)}(\boldsymbol{\theta}) := \{\mathrm{P}_{\boldsymbol{\theta};f}^{(n)}, f \in \mathcal{F}\})$ of $\mathcal{E}^{(n)}$, with maximal invariant the vector $\mathbf{R}^{(n)}(\boldsymbol{\theta}) := (R_1^{(n)}(\boldsymbol{\theta}), \dots, R_n^{(n)}(\boldsymbol{\theta}))'$, where $R_i^{(n)}(\boldsymbol{\theta})$ denotes the rank of the residual $Z_i^{(n)}(\boldsymbol{\theta})$ among $Z_1^{(n)}(\boldsymbol{\theta}), \dots, Z_n^{(n)}(\boldsymbol{\theta})$.

Similarly, let $\mathfrak{G}_+ := \{g \in \mathfrak{G} : g(-z) = -g(z)\}$ and denote by $\mathcal{G}_{\boldsymbol{\theta};+}^{(n)}$ the corresponding subgroup of $\mathcal{G}_{\boldsymbol{\theta}}^{(n)}$. This group (the group of *symmetric* order-preserving transformations of residuals) is a generating group for $\mathcal{E}_+^{(n)}(\boldsymbol{\theta}) := (\mathbb{R}^n, \mathcal{B}^n, \mathcal{P}_+^{(n)}(\boldsymbol{\theta}) := \{\mathrm{P}_{\boldsymbol{\theta};f}^{(n)}, f \in \mathcal{F}_+\})$, the submodel of $\mathcal{E}^{(n)}(\boldsymbol{\theta})$ that results from restricting to symmetric densities $f \in \mathcal{F}_+$. A maximal invariant here is the vector $\mathbf{R}_+^{(n)}(\boldsymbol{\theta}) := (s_1^{(n)}(\boldsymbol{\theta}) R_{+;1}^{(n)}(\boldsymbol{\theta}), \dots, s_n^{(n)}(\boldsymbol{\theta}) R_{+;n}^{(n)}(\boldsymbol{\theta}))'$, where $R_{+;i}^{(n)}(\boldsymbol{\theta})$ denotes the rank of the absolute value $|Z_i^{(n)}(\boldsymbol{\theta})|$ among $|Z_1^{(n)}(\boldsymbol{\theta})|, \dots, |Z_n^{(n)}(\boldsymbol{\theta})|$ and where $s_i^{(n)}(\boldsymbol{\theta})$ is the sign of $Z_i^{(n)}(\boldsymbol{\theta})$.

Turning to the model $\mathcal{E}_0^{(n)} := (\mathbb{R}^n, \mathcal{B}^n, \mathcal{P}_0^{(n)} := \{\mathrm{P}_{\boldsymbol{\theta};f}^{(n)}, \boldsymbol{\theta} \in \boldsymbol{\Theta}, f \in \mathcal{F}_0\})$ characterized by the residual function $\mathfrak{Z}_{\boldsymbol{\theta}}^{(n)}$ and the zero-median white noise concept $\mathcal{H}_{0;f}^{(n)}$, it is easy to see that a generating group for (with obvious notation) $\mathcal{E}_0^{(n)}(\boldsymbol{\theta})$ is obtained by considering the subgroup of $\mathcal{G}_{\boldsymbol{\theta}}^{(n)}$ that corresponds to $\mathfrak{G}_0 := \{g \in \mathfrak{G} : g(0) = 0\}$, with maximal invariant the vectors $\mathbf{s}^{(n)}(\boldsymbol{\theta}) := (s_1^{(n)}(\boldsymbol{\theta}), \dots, s_n^{(n)}(\boldsymbol{\theta}))'$ of residual signs and $\mathbf{R}^{(n)}(\boldsymbol{\theta})$ of residual ranks.

Provided that the parameter $\boldsymbol{\theta}$ contains a location or intercept component, and leaving aside the condition that residuals should have finite first-order moments, the model $\mathcal{E}_*^{(n)} := (\mathbb{R}^n, \mathcal{B}^n, \mathcal{P}_*^{(n)} := \{\mathrm{P}_{\boldsymbol{\theta};f}^{(n)}, \boldsymbol{\theta} \in \boldsymbol{\Theta}, f \in \mathcal{F}_*\})$, which is characterized by the same residual function $\mathfrak{Z}_{\boldsymbol{\theta}}^{(n)}$ as $\mathcal{E}_0^{(n)}$, but has zero-mean rather than zero-median white noise, coincides with $\mathcal{E}_0^{(n)}$. Both models indeed involve the same family of distributions $\mathcal{P}^{(n)}$ over $(\mathbb{R}^n, \mathcal{B}^n)$; they only differ in the way the nonparametric family $\mathcal{P}^{(n)}$ is split into a collection of parametric subfamilies $\mathcal{P}_f^{(n)} := \{\mathrm{P}_{f;\boldsymbol{\theta}}^{(n)}, \boldsymbol{\theta} \in \boldsymbol{\Theta}\}$ (hence, of course, in the way $\boldsymbol{\theta}$ is to be interpreted). Rather than two distinct models, $\mathcal{E}_0^{(n)}$ and $\mathcal{E}_*^{(n)}$ thus constitute two different parametrization of the same model, but the invariance structure underlying $\mathcal{E}_0^{(n)}$ is not present in $\mathcal{E}_*^{(n)}$. The median, in this respect, allows for a richer structure and, therefore, seems more appropriate than the mean as a location parameter.

2.3. *Group invariance and semiparametric efficiency.* The importance of considering maximal invariants—thus, signs and ranks in models with zero-median white noise—has been substantiated by Hallin and Werker [11]. Their paper showed that, in a very broad class of models, semiparametrically efficient inference procedures can be obtained by conditioning with respect to a maximal invariant $\sigma$-algebra.

More precisely, assume that the semiparametric family (2.1) is such that:

(i) For any fixed $f$, the parametric subfamily $\mathcal{P}_f^{(n)} := \{\mathrm{P}_{f;\boldsymbol{\theta}}^{(n)}, \boldsymbol{\theta} \in \boldsymbol{\Theta}\}$ is locally asymptotically normal (LAN), with central sequence $\boldsymbol{\Delta}_f^{(n)}(\boldsymbol{\theta})$.

(ii) For any fixed $\boldsymbol{\theta}$, the nonparametric subfamily $\mathcal{P}_{\boldsymbol{\theta}}^{(n)} := \{\mathrm{P}_{f;\boldsymbol{\theta}}^{(n)}, f \in \mathcal{F}\}$ is generated by a group of transformations with maximal invariant $\mathbf{W}^{(n)}(\boldsymbol{\theta})$.

Then, under very general conditions, semiparametrically efficient inference (testing, estimation, etc.) at $f$ can be based on the *semiparametrically efficient central sequence* $\mathrm{E}[\boldsymbol{\Delta}_f^{(n)}(\boldsymbol{\theta})|\mathbf{W}^{(n)}(\boldsymbol{\theta})]$, which, moreover, is distribution-free under $\mathcal{P}_{\boldsymbol{\theta}}^{(n)}$. Projecting onto maximal invariant $\sigma$-algebras (generated, in the context of Section 2.2, by the ranks, the signed ranks or the signs and ranks) thus yields (at given $f$) the same results as tangent space projections. In a companion paper [10], we specialize the Hallin and Werker [11] abstract results to obtain semiparametrically efficient inference in median regression and autoregressive models using the asymptotic representation results of the present paper for general sign-and-rank statistics.

Inference based on ranks and signed ranks has since long ago made its way to everyday practice and even to elementary textbooks. A pretty complete toolkit of rank-based methods is available for the analysis of linear models with independent observations (see [4, 18] for a systematic account and the state of the art in this context), as well as for the analysis of linear time series models (see [2, 5, 6, 7, 9]). It is somewhat surprising, therefore, that sign-and-rank statistics never have been considered so far in the vast literature devoted to that subject. The purpose of this paper is to fill this gap.

2.4. *Two simple examples.* Two examples are treated in some detail in Sections 3.4 (median regression) and 4.4 (median moving average), respectively.

Under the median-regression model, observations are of the form

$$Y_i = \theta_1 + \theta_2 c_i^{(n)} + \varepsilon_i, \qquad i = 1, \ldots, n, \tag{2.2}$$

where $\boldsymbol{\theta} := (\theta_1, \theta_2) \in \mathbb{R}^2$, the $c_i^{(n)}$'s are regression constants and the $\varepsilon_i$'s are independent and identically distributed (i.i.d.) with density $f$. Instead of the usual specification that $\mathrm{E}[\varepsilon_i] = 0$, however, we rather impose that the median of $\varepsilon_i$ is zero (i.e., $f \in \mathcal{F}_0$). Here, the residuals take the form $Z_i^{(n)}(\boldsymbol{\theta}) := Y_i - \theta_1 - \theta_2 c_i^{(n)}$. Under $\mathrm{P}_{f;\boldsymbol{\theta}}^{(n)}$, these residuals are i.i.d. with density $f \in \mathcal{F}_0$. Under fairly general conditions, this model, for fixed $f$ (with weak derivative $f'$), is LAN with central sequence

$$\boldsymbol{\Delta}_f^{(n)}(\boldsymbol{\theta}) := n^{-1/2} \sum_{i=1}^{n} \frac{-f'}{f}(Z_i^{(n)}(\boldsymbol{\theta})) \begin{pmatrix} 1 \\ c_i^{(n)} \end{pmatrix}. \tag{2.3}$$

In the first-order median moving average (MA) model, observations are generated by the MA equation

$$Y_t = \varepsilon_t + \theta \varepsilon_{t-1}, \qquad t = 1, \ldots, n, \tag{2.4}$$

with $\theta \in (-1, 1)$. Here again, we assume that the $\varepsilon_t$'s are independent and identically distributed with density $f$ and median zero. For simplicity, assume $\varepsilon_0 = 0$. The residuals are defined recursively as $Z_t^{(n)}(\boldsymbol{\theta}) := Y_t - \theta Z_{t-1}^{(n)}(\boldsymbol{\theta})$, with initial value $Z_0^{(n)}(\boldsymbol{\theta}) = 0$. Here again, for fixed $f$ (with weak derivative $f'$), LAN holds with central sequence

$$\boldsymbol{\Delta}_f^{(n)}(\boldsymbol{\theta}) := n^{-1/2} \sum_{t=1}^{n} \frac{-f'}{f}(Z_t^{(n)}(\boldsymbol{\theta})) Z_{t-1}^{(n)}(\boldsymbol{\theta}). \tag{2.5}$$

2.5. *Sign-and-rank statistics.* A *sign-and-rank statistic* is an $(\mathbf{s}^{(n)}, \mathbf{R}^{(n)})$-measurable statistic, where $\mathbf{s}^{(n)} = (s_1^{(n)}, \ldots, s_n^{(n)})'$ and $\mathbf{R}^{(n)} = (R_1^{(n)}, \ldots, R_n^{(n)})'$ are the vector of signs and the vector of ranks, respectively, associated with some $n$-dimensional random vector $\mathbf{Z}^{(n)}$. The objective of this paper is to introduce linear nonserial (Section 3) and linear serial (Section 4) sign-and-rank statistics, and to study their distributions under $\mathcal{H}_0^{(n)}$.

Denote by

$$N_-^{(n)} := \sum_{i=1}^{n} I[Z_i^{(n)} < 0] = \sum_{i=1}^{n} I[s_i^{(n)} = -1]$$

and by

$$N_+^{(n)} := \sum_{i=1}^{n} I[Z_i^{(n)} > 0] = \sum_{i=1}^{n} I[s_i^{(n)} = 1]$$

the numbers of negative and positive components in $\mathbf{Z}^{(n)}$ (in $\mathbf{s}^{(n)}$), respectively. Under $\mathcal{H}_0^{(n)}$, $N_+^{(n)}$ is binomial $\mathrm{Bin}(n, 1/2)$. Letting $\mathbf{N}^{(n)} := (N_-^{(n)}, N_+^{(n)})$, note that $\sigma(\mathbf{N}^{(n)}) = \sigma(N_-^{(n)}) = \sigma(N_+^{(n)})$, because $N_+^{(n)} = n - N_-^{(n)}$ with probability 1. Since $s_i^{(n)} = I[Z_i^{(n)} > 0] - I[Z_i^{(n)} < 0] = I[R_i^{(n)} > n - N_+^{(n)}] - I[R_i^{(n)} \leq N_-^{(n)}]$ for all $i = 1, \ldots, n$, the couple $(\mathbf{N}^{(n)}, \mathbf{R}^{(n)})$ is maximal invariant for $\mathcal{H}_0^{(n)}$.

Defining the sets

$$\mathcal{N}_-^{(n)} := \{i \in \{1,\ldots,n\} : s_i^{(n)} = -1\} = \{i_1^- < \cdots < i^-_{N_-^{(n)}}\}$$

and

$$\mathcal{N}_+^{(n)} := \{i \in \{1,\ldots,n\} : s_i^{(n)} = 1\} = \{i_1^+ < \cdots < i^+_{N_+^{(n)}}\},$$

the distribution of $(\mathbf{s}^{(n)}, \mathbf{R}^{(n)})$ under $\mathcal{H}_0^{(n)}$ is conveniently characterized as follows: The marginal distribution of $\mathbf{s}^{(n)}$ is uniform over the $2^n$ elements of $\{-1,1\}^n$ and the conditional distribution of $\mathbf{R}^{(n)}$ given $\mathbf{s}^{(n)}$ is such that $(R^{(n)}_{i_1^-}, R^{(n)}_{i_2^-}, \ldots, R^{(n)}_{i^-_{N_-^{(n)}}}; R^{(n)}_{i_1^+}, R^{(n)}_{i_2^+}, \ldots, R^{(n)}_{i^+_{N_+^{(n)}}})$ is (conditionally) uniformly distributed over the $(N_-^{(n)}!)(N_+^{(n)}!)$ possible combinations of a permutation of $\{1,\ldots,N_-^{(n)}\}$ with a permutation of $\{(n - N_+^{(n)}) + 1, \ldots, n\}$.

Let us finally denote by $\mathbf{Z}_{(\cdot)-}^{(N_-^{(n)})}$ and $\mathbf{Z}_{(\cdot)+}^{(N_+^{(n)})}$ the vectors of order statistics associated with the negative and positive elements of $\mathbf{Z}^{(n)}$, respectively. These two vectors—the first one of length $N_-^{(n)}$ and the second one of length $N_+^{(n)}$—constitute a natural (random) decomposition of the vector of order statistics $\mathbf{Z}_{(\cdot)}^{(n)}$ associated with $\mathbf{Z}^{(n)}$.

## 3. Nonserial linear sign-and-rank statistics.

3.1. *Definition and conditional asymptotic representation.* A *linear nonserial sign-and-rank statistic* is a statistic of the form

$$S_{\mathbf{c}}^{(n)} := \frac{1}{n}\sum_{i=1}^{n} c_i^{(n)} a^{(n)}(\mathbf{N}^{(n)}; R_i^{(n)}), \tag{3.1}$$

where $a^{(n)}(\cdot;\cdot)$ is a real-valued *score function* defined over $\{((\nu,\eta); i) : \nu, \eta \in \{0,1,\ldots,n\},\ \eta \le n-\nu,\ i \in \{1,\ldots,n\}\}$; note that each summand in (3.1) is allowed to depend on the sign $s_i^{(n)}$ of $Z_i^{(n)}$, but also, via $\mathbf{N}^{(n)}$, on the other signs, but not on the other ranks. As usual, the $c_i^{(n)}$'s $(i = 1,\ldots,n)$ denote nonrandom *regression constants*.

The exact mean $\mathrm{E}[S_{\mathbf{c}}^{(n)}]$ and the exact variance $\mathrm{Var}[S_{\mathbf{c}}^{(n)}]$ of $S_{\mathbf{c}}^{(n)}$ under $\mathcal{H}_0^{(n)}$ are easily obtained from elementary combinatorial arguments: Letting $\bar{c}^{(n)} := n^{-1}\sum_{i=1}^{n} c_i^{(n)}$, we obtain

$$\mathrm{E}[S_{\mathbf{c}}^{(n)}] = (n2^n)^{-1}\bar{c}^{(n)} \sum_{j=1}^{n}\sum_{\nu=0}^{n} \binom{n}{\nu} a^{(n)}((\nu, n-\nu); j)$$

and

$$\mathrm{Var}[S_{\mathbf{c}}^{(n)}] = \frac{1}{n(n-1)2^n} \sum_{i=1}^{n} (c_i^{(n)} - \bar{c}^{(n)})^2 \times \sum_{\nu=0}^{n} \binom{n}{\nu} \left\{ \sum_{i=1}^{n} [a^{(n)}((\nu, n-\nu); i)]^2 - \frac{1}{n} \left[ \sum_{i=1}^{n} a^{(n)}((\nu, n-\nu); i) \right]^2 \right\},$$

respectively.

If asymptotic results are to be obtained, some stability of the scores $a^{(n)}$ is required as $n$ increases. We therefore assume the existence of a *score-generating function*. A function $\varphi : (0,1) \to \mathbb{R}$ is called a *score-generating function* for the score function $a^{(n)}$ if

$$\mathrm{E}[\{a^{(n)}(\mathbf{N}^{(n)}; R_1^{(n)}) - \varphi(F(Z_1^{(n)}))\}^2 | \mathbf{Z}_{(\cdot)}^{(n)}] = o_{\mathrm{P}}(1) \tag{3.2}$$

under $\mathcal{H}_{0;f}^{(n)}$, as $n \to \infty$. Here $F$ denotes the distribution function associated with density $f$. Note that, by the rule of iterated expectations and the fact that $\mathbf{N}^{(n)} = (N_-^{(n)}, N_+^{(n)})$ is measurable with respect to $\mathbf{Z}_{(\cdot)}^{(n)}$, a sufficient condition for (3.2) to hold is

$$\mathrm{E}[\{a^{(n)}(\mathbf{N}^{(n)}; R_1^{(n)}) - \varphi(F(Z_1^{(n)}))\}^2 | \mathbf{N}^{(n)}] = o_{\mathrm{P}}(1) \tag{3.3}$$

under $\mathcal{H}_{0;f}^{(n)}$, as $n \to \infty$.

No asymptotic results for $S_{\mathbf{c}}^{(n)}$ can be obtained without some assumptions on the asymptotic behavior of regression constants $c_i^{(n)}, i = 1, \ldots, n$. We assume that the classical *Noether condition* holds:

(N) The constants $c_i^{(n)}$, $i = 1, \ldots, n$, are not all equal and

$$\lim_{n \to \infty} \frac{\max_{1 \le i \le n} (c_i^{(n)} - \bar{c}^{(n)})^2}{\sum_{j=1}^{n} (c_j^{(n)} - \bar{c}^{(n)})^2} = 0.$$

We may now state a first asymptotic representation and asymptotic normality result. This result, however, is a *conditional* one in the sense that the centering in (3.4) and (3.5) below is a conditional centering. Since, conditionally on the signs, the sign-and-rank statistic (3.1) reduces to a purely rank-based statistic, this conditional representation result follows from classical results on linear rank statistics and merely serves as an intermediate step in the derivation of the main result (of an unconditional nature) in Section 3.3. Contrary to the unconditional one, which requires *exact* or *approximate* scores, the conditional result holds for any scores that satisfy (3.2).

LEMMA 3.1. *Let* $\varphi:(0,1)\to\mathbb{R}$ *be a nonconstant square-integrable score-generating function for* $a^{(n)}$ *and let the regression constants* $c_i^{(n)}$ $(i=1,\dots,n)$ *satisfy the Noether condition* (N). *Assume moreover that* $\sum_{i=1}^n (c_i^{(n)}-\bar{c}^{(n)})^2 = O(n)$, *as* $n\to\infty$. *Then:*

(i) (Asymptotic representation) *under* $\mathcal{H}_{0;f}^{(n)}$, *as* $n\to\infty$,

$$S_{\mathbf{c}}^{(n)} - \mathrm{E}[S_{\mathbf{c}}^{(n)}|\mathbf{N}^{(n)}] = T_{\varphi;f}^{(n)} - \mathrm{E}[T_{\varphi;f}^{(n)}|\mathbf{Z}_{(\cdot)}^{(n)}] + o_{\mathrm{P}}(1/\sqrt{n}\,), \tag{3.4}$$

*where* $T_{\varphi;f}^{(n)} := \frac{1}{n}\sum_{i=1}^n c_i^{(n)}\varphi(F(Z_i^{(n)}))$ *(F stands for the distribution function associated with f);*

(ii) (Asymptotic normality) *under* $\mathcal{H}_0^{(n)}$, *as* $n\to\infty$,

$$\sqrt{n}(S_{\mathbf{c}}^{(n)} - \mathrm{E}[S_{\mathbf{c}}^{(n)}|\mathbf{N}^{(n)}])\Big/\sqrt{\frac{1}{n}\sum_{i=1}^n (c_i^{(n)}-\bar{c}^{(n)})^2} \xrightarrow{\mathcal{L}} \mathcal{N}(0,\sigma_\varphi^2), \tag{3.5}$$

*where* $0<\sigma_\varphi^2 := \int_0^1 \varphi^2(u)\,du - (\int_0^1\varphi(u)\,du)^2 < \infty$.

Observe that, under $\mathcal{H}_0^{(n)}$,

$$\begin{aligned}
\mathrm{E}[S_{\mathbf{c}}^{(n)}|\mathbf{N}^{(n)}] &= \frac{1}{n}\sum_{i=1}^n c_i^{(n)}\mathrm{E}[\mathrm{E}[a^{(n)}(\mathbf{N}^{(n)};R_i^{(n)})|\mathbf{s}^{(n)}]|\mathbf{N}^{(n)}] \\
&= \frac{1}{n}\sum_{i=1}^n c_i^{(n)}\Bigg\{\mathrm{P}[s_i^{(n)}=-1|\mathbf{N}^{(n)}]\frac{1}{N_-^{(n)}}\sum_{j=1}^{N_-^{(n)}} a^{(n)}(\mathbf{N}^{(n)};j) \\
&\qquad + \mathrm{P}[s_i^{(n)}=1|\mathbf{N}^{(n)}]\frac{1}{N_+^{(n)}}\sum_{j=(n-N_+^{(n)})+1}^{n} a^{(n)}(\mathbf{N}^{(n)};j)\Bigg\} \\
&= \bar{c}^{(n)}\left(\frac{1}{n}\sum_{j=1}^n a^{(n)}(\mathbf{N}^{(n)};j)\right) \\
&= \bar{c}^{(n)}\left(\frac{1}{n}\sum_{i=1}^n a^{(n)}(\mathbf{N}^{(n)};R_i^{(n)})\right)
\end{aligned}$$

and

$$\begin{aligned}
\mathrm{E}[T_{\varphi;f}^{(n)}|\mathbf{Z}_{(\cdot)}^{(n)}] &= \frac{1}{n}\sum_{i=1}^n c_i^{(n)}\mathrm{E}[\varphi(F(Z_i^{(n)}))|\mathbf{Z}_{(\cdot)}^{(n)}] \\
&= \bar{c}^{(n)}\left(\frac{1}{n}\sum_{i=1}^n \varphi(F(Z_i^{(n)}))\right).
\end{aligned} \tag{3.6}$$

Hence, part (i) of Lemma 3.1 actually states that

$$\begin{aligned}(3.7)\qquad &\frac{1}{n}\sum_{i=1}^{n}(c_i^{(n)}-\bar{c}^{(n)})a^{(n)}(\mathbf{N}^{(n)};R_i^{(n)})\\ &\quad=\frac{1}{n}\sum_{i=1}^{n}(c_i^{(n)}-\bar{c}^{(n)})\varphi(F(Z_i^{(n)}))+o_{\mathrm{P}}(1/\sqrt{n}\,),\end{aligned}$$

under $\mathcal{H}_{0;f}^{(n)}$, as $n\to\infty$. Note that the expression on the right-hand side of (3.7) coincides with the asymptotic representation of the purely rank-based statistic $\frac{1}{n}\sum_{i=1}^{n}(c_i^{(n)}-\bar{c}^{(n)})a_\varphi^{(n)}(R_i^{(n)})$, where $a_\varphi^{(n)}(R_i^{(n)})$ are, for instance, the traditional exact scores $\mathrm{E}[\varphi(F(Z_i^{(n)}))|R_i^{(n)}]$ associated with the score-generating function $\varphi$. The sign-and-rank statistic $S_{\mathbf{c}}^{(n)}$ thus asymptotically decomposes into two parts; one of them (namely, $S_{\mathbf{c}}^{(n)}-\mathrm{E}[S_{\mathbf{c}}^{(n)}|\mathbf{N}^{(n)}]$) asymptotically does not depend on $\mathbf{N}^{(n)}$ and represents the contribution of the ranks, while the second one ($\mathrm{E}[S_{\mathbf{c}}^{(n)}|\mathbf{N}^{(n)}]-\mathrm{E}[S_{\mathbf{c}}^{(n)}]$) constitutes the contribution of the signs. Moreover, the ranks and $\mathbf{N}^{(n)}$ being mutually independent, these two quantities are orthogonal to each other and contribute additively to the unconditional asymptotic variance (see the proof of Proposition 3.2 below).

PROOF OF LEMMA 3.1. Since the ranks $\mathbf{R}^{(n)}$ and $\mathbf{N}^{(n)}$ are mutually independent under $\mathcal{H}_{0;f}^{(n)}$, part (i) of the lemma follows from classical asymptotic representation results for linear rank statistics; see [3], page 61. The proof of part (ii) of the lemma, in view of (3.4), simply consists in checking that $\sqrt{n}(T_{\varphi;f}^{(n)}-\mathrm{E}[T_{\varphi;f}^{(n)}|\mathbf{Z}_{(\cdot)}^{(n)}])$ satisfies the traditional Lindeberg condition. □

3.2. *Exact and approximate scores.* Following the classical literature on ranks, we consider in the present paper sign-and-rank statistics based on either *exact* or *approximate* scores.

Let $U_1^{(n)},\dots,U_n^{(n)}$ be an $n$-tuple of i.i.d. random variables uniformly distributed over $(0,1)$. Define $s_{U_i}^{(n)}:=I[U_i^{(n)}>1/2]-I[U_i^{(n)}<1/2]$, $N_{\mathbf{U};-}^{(n)}:=\sum_{i=1}^{n}I[U_i^{(n)}<1/2]$ and $N_{\mathbf{U};+}^{(n)}:=\sum_{i=1}^{n}I[U_i^{(n)}>1/2]$. Denote by $R_{U_i}^{(n)}$ the rank of $U_i^{(n)}$ among $U_1^{(n)},\dots,U_n^{(n)}$, by $U_{(i)-}^{(\nu)}$ $(i=1,\dots,\nu)$ the $i$th-order statistic associated with a sample of $\nu$ i.i.d. random variables uniformly distributed over $(0,1/2)$ and by $U_{(i)+}^{(\nu)}$ $(i=1,\dots,\nu)$ the $i$th-order statistic associated with a sample of $\nu$ i.i.d. random variables uniformly distributed over $(1/2,1)$. Note that the conditional distribution of $U_i^{(n)}$ given the event $s_{U_i}^{(n)}=-1$ (resp. $s_{U_i}^{(n)}=1$) is uniform over $(0,1/2)$ [resp. $(1/2,1)$]. The linear nonserial sign-and-rank statistics constructed from the *exact* and *approximate* scores associated with $\varphi$ are defined by

$$\begin{aligned} S^{(n)}_{\mathbf{c};\varphi;\mathrm{ex/appr}} &:= \frac{1}{n}\sum_{i=1}^{n} c_i^{(n)} a^{(n)}_{\varphi;\mathrm{ex/appr}}(\mathbf{N}^{(n)}; R_i^{(n)}) \\ &:= \frac{1}{n}\sum_{i=1}^{n} c_i^{(n)}\{I[s_i^{(n)}=-1]a^{(n)}_{\varphi;-;\mathrm{ex/appr}}(N_-^{(n)}; R_i^{(n)}) \\ &\qquad + I[s_i^{(n)}=1]a^{(n)}_{\varphi;+;\mathrm{ex/appr}}(N_+^{(n)}; R_i^{(n)}-(n-N_+^{(n)}))\}, \end{aligned} \tag{3.8}$$

where the score functions $a^{(n)}_{\varphi;-;\mathrm{ex}}$, $a^{(n)}_{\varphi;-;\mathrm{appr}}$, $a^{(n)}_{\varphi;+;\mathrm{ex}}$ and $a^{(n)}_{\varphi;+;\mathrm{appr}}$, all defined on the set $\{(\nu;i);\nu,i\in\{1,\dots,n\}$ with $i\le\nu\}$, are given by

$$a^{(n)}_{\varphi;-;\mathrm{ex}}(\nu;i) := \mathrm{E}[\varphi(U_1^{(n)})|N^{(n)}_{\mathbf{U};-}=\nu, R^{(n)}_{U_1}=i] = \mathrm{E}[\varphi(U^{(\nu)}_{(i)-})], \tag{3.9}$$

$$a^{(n)}_{\varphi;-;\mathrm{appr}}(\nu;i) := \varphi(\mathrm{E}[U^{(\nu)}_{(i)-}]) = \varphi\left(\frac{i}{2(\nu+1)}\right), \tag{3.10}$$

$$\begin{aligned} a^{(n)}_{\varphi;+;\mathrm{ex}}(\nu;i) &:= \mathrm{E}[\varphi(U_1^{(n)})|N^{(n)}_{\mathbf{U};+}=\nu, R^{(n)}_{U_1}=(n-\nu)+i] \\ &= \mathrm{E}[\varphi(U^{(\nu)}_{(i)+})] \end{aligned} \tag{3.11}$$

and

$$a^{(n)}_{\varphi;+;\mathrm{appr}}(\nu;i) := \varphi(\mathrm{E}[U^{(\nu)}_{(i)+}]) = \varphi\left(\frac{1}{2}+\frac{i}{2(\nu+1)}\right) \tag{3.12}$$

Observe that, under $\mathcal{H}^{(n)}_{0;f}$, $S^{(n)}_{\mathbf{c};\varphi;\mathrm{ex}} = \mathrm{E}[T^{(n)}_{\varphi;f}|\mathbf{N}^{(n)},\mathbf{R}^{(n)}] = \mathrm{E}[T^{(n)}_{\varphi;f}|\mathbf{s}^{(n)},\mathbf{R}^{(n)}]$.

We then have the following proposition.

PROPOSITION 3.1. *Let $\varphi:(0,1)\to\mathbb{R}$ be a nonconstant square-integrable function. Then $\varphi$ is a score-generating function for $a^{(n)}_{\varphi;\mathrm{ex}}$. If, moreover, $\varphi$ is the difference of two nondecreasing square-integrable functions, then $\varphi$ is also a score-generating function for $a^{(n)}_{\varphi;\mathrm{appr}}$.*

PROOF. Let us first consider the *exact* scores defined by relationships (3.8), (3.9) and (3.11), and let us show that, under $\mathcal{H}^{(n)}_{0;f}$,

$$\mathrm{E}[\{a^{(n)}_{\varphi;\mathrm{ex}}(\mathbf{N}^{(n)}; R_1^{(n)}) - \varphi(F(Z_1^{(n)}))\}^2|\mathbf{N}^{(n)}] = o_{\mathrm{P}}(1) \tag{3.13}$$

as $n\to\infty$. By the definition of $a^{(n)}_{\varphi;-;\mathrm{ex}}$, we only need to show that

$$\begin{aligned} &\mathrm{E}[\{\mathrm{E}[\varphi(F(Z_1^{(n)}))|s_1^{(n)}=-1, N_-^{(n)}, R_1^{(n)}] - \varphi(F(Z_1^{(n)}))\}^2|N_-^{(n)}, s_1^{(n)}=-1] \\ &\quad = o_{\mathrm{P}}(1), \end{aligned}$$

under $\mathcal{H}_{0;f}^{(n)}$, as $n \to \infty$. Since $F(Z_1^{(n)})$ is, under $\mathcal{H}_{0;f}^{(n)}$ and conditionally on $s_1^{(n)} = -1$, uniform over the interval $(0, 1/2)$, this readily follows from a slight generalization of Theorem V.1.4.a in [3], page 157.

Let us now consider the *approximate* scores defined by (3.8), (3.10) and (3.12). Clearly, (3.3) holds for $a_{\varphi;\mathrm{appr}}^{(n)}$ if, under $\mathcal{H}_{0;f}^{(n)}$,

$$\mathrm{E}[\{a_{\varphi;-;\mathrm{appr}}^{(n)}(N_-^{(n)}; R_1^{(n)}) - \varphi(F(Z_1^{(n)}))\}^2 | N_-^{(n)}, s_1^{(n)} = -1]$$

and

$$\mathrm{E}[\{a_{\varphi;+;\mathrm{appr}}^{(n)}(N_+^{(n)}; R_1^{(n)} - (n - N_+^{(n)})) - \varphi(F(Z_1^{(n)}))\}^2 | N_+^{(n)}, s_1^{(n)} = 1]$$

are $o_{\mathrm{P}}(1)$ as $n \to \infty$. We have

$$\begin{aligned}
&\mathrm{E}[\{a_{\varphi;-;\mathrm{appr}}^{(n)}(N_-^{(n)}; R_1^{(n)}) - \varphi(F(Z_1^{(n)}))\}^2 | N_-^{(n)}, s_1^{(n)} = -1] \\
&\quad = \mathrm{E}[\{(a_{\varphi;-;\mathrm{appr}}^{(n)}(N_-^{(n)}; R_1^{(n)}) - a_{\varphi;-;\mathrm{ex}}^{(n)}(N_-^{(n)}; R_1^{(n)})) \\
&\qquad\qquad + (a_{\varphi;-;\mathrm{ex}}^{(n)}(N_-^{(n)}; R_1^{(n)}) - \varphi(F(Z_1^{(n)})))\}^2 | N_-^{(n)}, s_1^{(n)} = -1] \\
&\quad \le 2\mathrm{E}[\{a_{\varphi;-;\mathrm{appr}}^{(n)}(N_-^{(n)}; R_1^{(n)}) - a_{\varphi;-;\mathrm{ex}}^{(n)}(N_-^{(n)}; R_1^{(n)})\}^2 | N_-^{(n)}, s_1^{(n)} = -1] \\
&\qquad + 2\mathrm{E}[\{a_{\varphi;-;\mathrm{ex}}^{(n)}(N_-^{(n)}; R_1^{(n)}) - \varphi(F(Z_1^{(n)}))\}^2 | N_-^{(n)}, s_1^{(n)} = -1].
\end{aligned}$$

In view of the result for exact scores, we just consider the second term. Denoting by $\lfloor x \rfloor$ the integer part of $x$ ($x \in \mathbb{R}^+$), we may write

$$\begin{aligned}
&\mathrm{E}[\{a_{\varphi;-;\mathrm{appr}}^{(n)}(N_-^{(n)}; R_1^{(n)}) - a_{\varphi;-;\mathrm{ex}}^{(n)}(N_-^{(n)}; R_1^{(n)})\}^2 | N_-^{(n)}, s_1^{(n)} = -1] \\
&\quad = \frac{1}{N_-^{(n)}} \sum_{i=1}^{N_-^{(n)}} \{a_{\varphi;-;\mathrm{appr}}^{(n)}(N_-^{(n)}; i) - a_{\varphi;-;\mathrm{ex}}^{(n)}(N_-^{(n)}; i)\}^2 \\
&\quad = \int_0^1 \{(a_{\varphi;-;\mathrm{appr}}^{(n)}(N_-^{(n)}; 1 + \lfloor N_-^{(n)} u \rfloor) - \varphi(u/2)) \\
&\qquad\qquad + (\varphi(u/2) - a_{\varphi;-;\mathrm{ex}}^{(n)}(N_-^{(n)}; 1 + \lfloor N_-^{(n)} u \rfloor))\}^2 \, du \\
&\quad \le 2 \int_0^1 \{a_{\varphi;-;\mathrm{appr}}^{(n)}(N_-^{(n)}; 1 + \lfloor N_-^{(n)} u \rfloor) - \varphi(u/2)\}^2 \, du \\
&\qquad + 2 \int_0^1 \{a_{\varphi;-;\mathrm{ex}}^{(n)}(N_-^{(n)}; 1 + \lfloor N_-^{(n)} u \rfloor) - \varphi(u/2)\}^2 \, du.
\end{aligned}$$

That this latter quantity is $o_{\mathrm{P}}(1)$ follows from an obvious adaptation of Lemma V.1.6.a and Theorem V.1.4.b in [3], pages 164 and 158, respectively. □

3.3. *Asymptotic representation and asymptotic normality.* We now can state, for the nonserial case, the main result of this paper.

PROPOSITION 3.2. *Let $\varphi:(0,1)\to\mathbb{R}$ be a nonconstant square-integrable score-generating function for $S^{(n)}_{\mathbf{c};\varphi;\mathrm{ex/appr}}$ and let the regression constants $c_i^{(n)}$ $(i=1,\ldots,n)$ satisfy the Noether condition* (N). *Whenever approximate scores are considered, assume that $\varphi$ is the difference of two nondecreasing square-integrable functions. Assume, moreover, that $\bar{c}^{(n)}=O(1)$ and $\sum_{i=1}^n (c_i^{(n)}-\bar{c}^{(n)})^2=O(n)$ as $n\to\infty$. Let $\mu_\varphi^- := \int_0^{1/2}\varphi(u)\,du$, $\mu_\varphi^+ := \int_{1/2}^1 \varphi(u)\,du$ and $\mu_\varphi := \int_0^1\varphi(u)\,du$. Then, writing $S^{(n)}_{\mathbf{c}}$ for either $S^{(n)}_{\mathbf{c};\varphi;\mathrm{ex}}$ or $S^{(n)}_{\mathbf{c};\varphi;\mathrm{appr}}$:*

(i) (Asymptotic representation) *under $\mathcal{H}^{(n)}_{0;f}$, as $n\to\infty$,*

$$\begin{aligned} S^{(n)}_{\mathbf{c}} - \mathrm{E}[S^{(n)}_{\mathbf{c}}] &= \frac{1}{n}\sum_{i=1}^n (c_i^{(n)}-\bar{c}^{(n)})\varphi(F(Z_i^{(n)})) \\ &\quad + \bar{c}^{(n)}\left\{2\frac{N_-^{(n)}}{n}\mu_\varphi^- + 2\frac{N_+^{(n)}}{n}\mu_\varphi^+ - \mu_\varphi\right\} + o_{\mathrm{P}}(1/\sqrt{n}); \end{aligned} \tag{3.14}$$

(ii) (Asymptotic normality) *under $\mathcal{H}^{(n)}_0$, as $n\to\infty$,*

$$\begin{aligned} &\sqrt{n}(S^{(n)}_{\mathbf{c}} - \mathrm{E}[S^{(n)}_{\mathbf{c}}])\Big/\sqrt{\frac{\sigma_\varphi^2}{n}\sum_{i=1}^n (c_i^{(n)}-\bar{c}^{(n)})^2 + [\bar{c}^{(n)}(\mu_\varphi^- - \mu_\varphi^+)]^2} \\ &\quad \xrightarrow{\mathcal{L}} \mathcal{N}(0,1). \end{aligned} \tag{3.15}$$

Note that, in case $\varphi$ is skew-symmetric with respect to $1/2$ [i.e., $\varphi(u)=-\varphi(1-u)$], we have $\mu_\varphi^- = -\mu_\varphi^+$ and $\mu_\varphi=0$. Straightforward calculation yields $\bar{c}^{(n)}\{2\frac{N_-^{(n)}}{n}\mu_\varphi^- + 2\frac{N_+^{(n)}}{n}\mu_\varphi^+ - \mu_\varphi\} = \bar{c}^{(n)}\mu_\varphi^-(1-2\frac{N_-^{(n)}}{n})$. The conditional (3.4) and unconditional (3.14) asymptotic representations thus coincide and reduce to Hájek's traditional one for linear rank statistics, as soon as $\bar{c}^{(n)}=o(1)$ (examples of skew-symmetric score functions are the location scores $\varphi_f := -f'/f$ of a symmetric distribution with density $f$).

PROOF OF PROPOSITION 3.2. (i) We first establish (3.14) for exact scores. From (3.4) and (3.6), we have

$$\begin{aligned} S^{(n)}_{\mathbf{c};\varphi;\mathrm{ex}} - \mathrm{E}[S^{(n)}_{\mathbf{c};\varphi;\mathrm{ex}}] &= \frac{1}{n}\sum_{i=1}^n (c_i^{(n)}-\bar{c}^{(n)})\varphi(F(Z_i^{(n)})) \\ &\quad + \mathrm{E}[S^{(n)}_{\mathbf{c};\varphi;\mathrm{ex}}|\mathbf{N}^{(n)}] - \mathrm{E}[S^{(n)}_{\mathbf{c};\varphi;\mathrm{ex}}] + o_{\mathrm{P}}(1/\sqrt{n}). \end{aligned} \tag{3.16}$$

Since

$$\begin{aligned}\mathrm{E}[S^{(n)}_{\mathbf{c};\varphi;\mathrm{ex}}|\mathbf{N}^{(n)}] &= \mathrm{E}[\mathrm{E}[T^{(n)}_{\varphi;f}|\mathbf{N}^{(n)},\mathbf{R}^{(n)}]|\mathbf{N}^{(n)}] = \mathrm{E}[T^{(n)}_{\varphi;f}|\mathbf{N}^{(n)}]\\ &= \mathrm{E}[\mathrm{E}[T^{(n)}_{\varphi;f}|\mathbf{s}^{(n)}]|\mathbf{N}^{(n)}]\\ &= \mathrm{E}\left[\frac{1}{n}\sum_{i=1}^{n} c^{(n)}_i \mathrm{E}[\varphi(F(Z^{(n)}_i))|s^{(n)}_i]|\mathbf{N}^{(n)}\right],\end{aligned}$$

where

$$\begin{aligned}\mathrm{E}[\varphi(F(Z^{(n)}_i))|s^{(n)}_i] &= I[s^{(n)}_i = -1]\int_0^{1/2}\varphi(u)2\,du\\ &\quad + I[s^{(n)}_i = 1]\int_{1/2}^{1}\varphi(u)2\,du\\ &= 2I[s^{(n)}_i = -1]\mu^-_\varphi + 2I[s^{(n)}_i = 1]\mu^+_\varphi,\end{aligned}$$

it follows that

$$\begin{aligned}\mathrm{E}[S^{(n)}_{\mathbf{c};\varphi;\mathrm{ex}}|\mathbf{N}^{(n)}] &= \frac{2}{n}\sum_{i=1}^{n} c^{(n)}_i \mathrm{E}[I[s^{(n)}_i = -1]\mu^-_\varphi\\ &\qquad + I[s^{(n)}_i = 1]\mu^+_\varphi|\mathbf{N}^{(n)}]\\ &= 2\bar{c}^{(n)}\left(\frac{N^{(n)}_-}{n}\mu^-_\varphi + \frac{N^{(n)}_+}{n}\mu^+_\varphi\right)\end{aligned}\tag{3.17}$$

and

$$\mathrm{E}[S^{(n)}_{\mathbf{c};\varphi;\mathrm{ex}}|\mathbf{N}^{(n)}] - \mathrm{E}[S^{(n)}_{\mathbf{c};\varphi;\mathrm{ex}}] = \bar{c}^{(n)}\left(2\frac{N^{(n)}_-}{n}\mu^-_\varphi + 2\frac{N^{(n)}_+}{n}\mu^+_\varphi - \mu_\varphi\right),\tag{3.18}$$

which, along with (3.16), establish (3.14) for exact scores.

Turning to approximate scores, we can assume, without loss of generality, that $\varphi$ is nondecreasing. Since (3.16) also holds if approximate scores are substituted for the exact ones, it is sufficient, so that (3.14) holds for approximate scores, to show that the difference

$$\begin{aligned}E^{(n)} &:= \{\mathrm{E}[S^{(n)}_{\mathbf{c};\varphi;\mathrm{appr}}|\mathbf{N}^{(n)}] - \mathrm{E}[S^{(n)}_{\mathbf{c};\varphi;\mathrm{appr}}]\}\\ &\qquad - \{\mathrm{E}[S^{(n)}_{\mathbf{c};\varphi;\mathrm{ex}}|\mathbf{N}^{(n)}] - \mathrm{E}[S^{(n)}_{\mathbf{c};\varphi;\mathrm{ex}}]\}\end{aligned}\tag{3.19}$$

is $o_{\mathrm{P}}(1/\sqrt{n})$. Note that

$$\begin{aligned}&\mathrm{E}[S^{(n)}_{\mathbf{c};\varphi;\mathrm{appr}}|\mathbf{N}^{(n)}]\\ &\quad = \bar{c}^{(n)}\frac{1}{n}\left\{\sum_{j=1}^{N^{(n)}_-}\varphi\left(\frac{j}{2(N^{(n)}_- + 1)}\right) + \sum_{j=1}^{N^{(n)}_+}\varphi\left(\frac{1}{2} + \frac{j}{2(N^{(n)}_+ + 1)}\right)\right\}\end{aligned}\tag{3.20}$$

$$= \bar{c}^{(n)}\left\{2\frac{N_-^{(n)}}{n}D^-_{N_-^{(n)}} + 2\frac{N_+^{(n)}}{n}D^+_{N_+^{(n)}}\right\},$$

where $D_m^- := \frac{1}{2m}\sum_{j=1}^m \varphi(\frac{j}{2(m+1)})$ and $D_m^+ := \frac{1}{2m}\sum_{j=1}^m \varphi(\frac{1}{2}+\frac{j}{2(m+1)})$ are Riemann sums for the integrals $\mu_\varphi^- := \int_0^{1/2}\varphi(u)\,du$ and $\mu_\varphi^+ := \int_{1/2}^1 \varphi(u)\,du$, respectively. Since $\varphi$ is square-integrable, any term in the Riemann sum $\frac{1}{2m}\sum_{j=1}^m \varphi^2(\frac{1}{2}+\frac{j}{2(m+1)})$ associated with $\int_{1/2}^1 \varphi^2(u)\,du$ is $o(1)$ as $m\to\infty$. This implies that $\frac{1}{2m}\varphi(\frac{1}{2}+\frac{m}{2(m+1)})$ is $o(1/\sqrt{m})$; hence, in view of the fact that $N_+^{(n)} = O_{\mathrm{P}}(n)$, this implies that $\frac{1}{2N_+^{(n)}}\varphi(\frac{1}{2}+\frac{N_+^{(n)}}{2(N_+^{(n)}+1)}) = o_{\mathrm{P}}(1/\sqrt{n}\,)$ as $n\to\infty$. The same reasoning shows that any finite sum of Riemann terms in $D^-_{N_-^{(n)}}$ or $D^+_{N_+^{(n)}}$ actually is $o_{\mathrm{P}}(1/\sqrt{n}\,)$ as $n\to\infty$.

Now, any Riemann sum $D_m^+$ for $\mu_\varphi^+$ satisfies, since $\varphi$ is nondecreasing, the double inequality $\underline{D}_m^+ \le D_m^+ \le \bar{D}_m^+$, where $\underline{D}_m^+ := \frac{1}{2m}\sum_{j=0}^{m-1}\varphi(\frac{1}{2}+\frac{j}{2(m+1)})$ and $\bar{D}_m^+ := \frac{1}{2m}\sum_{j=1}^m \varphi(\frac{1}{2}+\frac{j}{2(m+1)})$ are the upper and lower Darboux sums associated with $\int_{1/2}^1\varphi(u)\,du$. The difference $\bar{D}_m^+ - \underline{D}_m^+$ clearly is $\frac{1}{2m}(\varphi(\frac{1}{2}+\frac{m}{2(m+1)}) - \varphi(\frac{1}{2}))$, which is $o(1/\sqrt{m}\,)$ as $m\to\infty$. Hence, for any Riemann sum, $D_m^+ - \mu_\varphi^+$ is also $o(1/\sqrt{m}\,)$, so that $D^+_{N_+^{(n)}} - \mu_\varphi^+ = o_{\mathrm{P}}(1/\sqrt{n}\,)$ as $n\to\infty$.

Furthermore, since the sequence $D_m^+ - \mu_\varphi^+$ converges to zero, it is bounded, so that $D^+_{N_+^{(n)}} - \mu_\varphi^+$ is uniformly integrable and $\mathrm{E}[\frac{N_+^{(n)}}{n}D^+_{N_+^{(n)}} - \frac{1}{2}\mu_\varphi^+] = o(1/\sqrt{n}\,)$ as $n\to\infty$.

A similar reasoning of course holds for $D^-_{N_-^{(n)}}$ and $\mu_\varphi^-$. Going back to (3.20) and recalling that $\bar{c}^{(n)} = O(1)$, we thus obtain the desired result that $E^{(n)}$ is $o_{\mathrm{P}}(1/\sqrt{n}\,)$. This completes the proof of part (i) of the proposition.

(ii) As for asymptotic normality, elementary calculations yield

$$\sqrt{n}\bar{c}^{(n)}\left(2\frac{N_-^{(n)}}{n}\mu_\varphi^- + 2\frac{N_+^{(n)}}{n}\mu_\varphi^+ - \mu_\varphi\right)$$
$$= \bar{c}^{(n)}\left(2(\mu_\varphi^- - \mu_\varphi^+)\left(\frac{N_-^{(n)}}{n} - \frac{1}{2}\right)\Big/\sqrt{1/4n}\right)\sqrt{\frac{1}{4}},$$

which, since $(\frac{N_-^{(n)}}{n} - \frac{1}{2})/\sqrt{1/4n}$ is asymptotically standard normal, is also asymptotically normal with mean zero and asymptotic variance $[\bar{c}^{(n)}(\mu_\varphi^- - \mu_\varphi^+)]^2$. The remark (right after Lemma 3.1) on the orthogonality between the two parts of the asymptotic representation of $S_{\mathbf{c}}^{(n)}$ completes the proof. □

Test statistics related to "regression coefficients" naturally involve "regression constants" $c_i^{(n)}$ that are not all equal. Quite on the contrary, test statistics related to location and intercepts do not involve any constants—more precisely, they are still of the form $S_{\mathbf{c}}^{(n)}$, but with constants $c_i^{(n)}$ all equal to 1. Proposition 3.2, as it is stated, does not apply. However, going back to the proof, one easily checks that, letting $S_{\varphi;\mathrm{ex/appr}}^{(n)} := \frac{1}{n}\sum_{i=1}^n a_{\varphi;\mathrm{ex/appr}}^{(n)}(\mathbf{N}^{(n)}; R_i^{(n)})$, under the same assumptions on the scores $\varphi$,

$$\begin{aligned} &\sqrt{n}(S_{\varphi;\mathrm{ex/appr}}^{(n)} - \mathrm{E}[S_{\varphi;\mathrm{ex/appr}}^{(n)}]) \\ &\quad = 2\frac{N_-^{(n)}}{n}\mu_\varphi^- + 2\frac{N_+^{(n)}}{n}\mu_\varphi^+ - \mu_\varphi + o_{\mathrm{P}}(1) \\ &\quad \xrightarrow{\mathcal{L}} \mathcal{N}(0, (\mu_\varphi^- - \mu_\varphi^+)^2) \end{aligned} \tag{3.21}$$

under $\mathcal{H}_0^{(n)}$, as $n \to \infty$.

3.4. *Example: median regression.* The central sequence (2.3) takes the form $\mathbf{\Delta}_f^{(n)} = \mathbf{\Delta}_f^{(n)}(\boldsymbol{\theta}) = \sqrt{n}(T_{\varphi_f;f;1}^{(n)}, T_{\varphi_f;f;2}^{(n)})'$ with (using the notation of Section 3)

$$T_{\varphi_f;f;1}^{(n)} := \frac{1}{n}\sum_{i=1}^n \varphi_f(F(Z_i^{(n)})),$$

$$T_{\varphi_f;f;2}^{(n)} := \frac{1}{n}\sum_{i=1}^n c_i^{(n)}\varphi_f(F(Z_i^{(n)}))$$

and $\varphi_f(u) := \frac{-f'}{f}(F^{-1}(u))$, $u \in (0,1)$. Instead of an arbitrary score-generating function, we therefore focus on $\varphi_f$. Define

$$S_{\varphi_f;1;\mathrm{ex}}^{(n)} := \mathrm{E}[T_{\varphi_f;f;1}^{(n)} | \mathbf{N}^{(n)}, \mathbf{R}^{(n)}]$$

and

$$S_{\varphi_f;2;\mathrm{ex}}^{(n)} := \mathrm{E}[T_{\varphi_f;f;2}^{(n)} | \mathbf{N}^{(n)}, \mathbf{R}^{(n)}].$$

Straightforward calculations lead to

$$2\frac{N_-^{(n)}}{n}\mu_{\varphi_f}^- + 2\frac{N_+^{(n)}}{n}\mu_{\varphi_f}^+ - \mu_{\varphi_f} = 2f(0)\frac{N_+^{(n)} - N_-^{(n)}}{n},$$

so that Proposition 3.2 and (3.21) yield

$$\mathrm{E}[\mathbf{\Delta}_f^{(n)} | \mathbf{N}^{(n)}, \mathbf{R}^{(n)}] = \sqrt{n}\begin{pmatrix} S_{\varphi_f;1;\mathrm{ex}}^{(n)} \\ S_{\varphi_f;2;\mathrm{ex}}^{(n)} \end{pmatrix} = \mathbf{\Delta}_f^{(n)*} + o_{\mathrm{P}}(1)$$

under $\mathrm{P}^{(n)}_{f;\boldsymbol{\theta}}$, as $n \to \infty$, where

$$\boldsymbol{\Delta}^{(n)*}_{f} := \sqrt{n}\begin{pmatrix} 2f(0)\dfrac{N^{(n)}_{+} - N^{(n)}_{-}}{n} \\ \dfrac{1}{n}\sum_{i=1}^{n}(c^{(n)}_{i} - \bar{c}^{(n)})\varphi_f(F(Z^{(n)}_{i})) + \bar{c}^{(n)}2f(0)\dfrac{N^{(n)}_{+} - N^{(n)}_{-}}{n} \end{pmatrix}$$

is a version of the semiparametrically efficient central sequence for $\boldsymbol{\theta}$ in the semiparametric experiment $\mathcal{E}^{(n)}_{0}$. This latter statement can easily be checked using standard tangent space calculations. Similarly, in view of (3.5) and Proposition 3.1, the approximate score version of the same semiparametrically efficient central sequence is

$$\underset{\sim}{\Delta}{}^{(n)*}_{f} := \sqrt{n}\begin{pmatrix} 2f(0)\dfrac{N^{(n)}_{+} - N^{(n)}_{-}}{n} \\ \dfrac{1}{n}\sum_{i=1}^{n}(c^{(n)}_{i} - \bar{c}^{(n)})\varphi_f(\underset{\sim}{R}{}^{(n)}_{i}) + \bar{c}^{(n)}2f(0)\dfrac{N^{(n)}_{+} - N^{(n)}_{-}}{n} \end{pmatrix}$$

with, for $i = 1, \ldots, n$,

$$\begin{aligned} \underset{\sim}{R}{}^{(n)}_{i} &:= I[R^{(n)}_{i} \le N^{(n)}_{-}]\frac{R^{(n)}_{i}}{2(N^{(n)}_{-} + 1)} \\ &\quad + I[R^{(n)}_{i} > n - N^{(n)}_{+}]\left(\frac{1}{2} + \frac{R^{(n)}_{i} - (n - N^{(n)}_{+})}{2(N^{(n)}_{+} + 1)}\right). \end{aligned} \tag{3.22}$$

This central sequence, which is measurable with respect to the residual signs and ranks, can be used to perform semiparametrically efficient inference (tests, estimation, etc.); see, for example, Section 11.9 of [16]. For a full treatment of sign-and-rank-based versions of semiparametrically efficient central sequences in median restricted models, we refer to [10].

## 4. Serial linear sign-and-rank statistics.

4.1. *Definition and conditional asymptotic representation.* Nonserial sign-and-rank statistics, just as their traditional rank-based counterparts, are inefficient in the context of dependent observations: Only serial statistics can capture the effects of serial dependence. Define a linear serial sign-and-rank statistic of order $k$ ($k \in \{1, \ldots, n-1\}$) as a statistic of the form

$$S^{(n)}_{k} := \frac{1}{n-k}\sum_{t=k+1}^{n} a^{(n)}_{k}(\mathbf{N}^{(n)}; R^{(n)}_{t}, \ldots, R^{(n)}_{t-k}),$$

where $a_k^{(n)}(\cdot;\cdot,\dots,\cdot)$ is defined over the product of the set $\{(\nu,\eta);\nu,\eta\in\{0,1,\dots,n\},\eta\le n-\nu\}$ with the set of all $(k+1)$-tuples of distinct integers in $\{1,\dots,n\}$. The asymptotic mean and variance of $S_k^{(n)}$ are given in the subsequent Proposition 4.1.

Here also an asymptotic representation result is proved, establishing the asymptotic equivalence between $S_k^{(n)}$ and a "parametric" serial statistic $T_k^{(n)}$. The asymptotic normality of $T_k^{(n)}$ then entails that of $S_k^{(n)}$. A function $\varphi_k:(0,1)^{k+1}\to\mathbb{R}$ is a score-generating function for the serial score function $a_k^{(n)}$ if

$$\begin{aligned}&\mathrm{E}[\{a_k^{(n)}(\mathbf{N}^{(n)};R_{k+1}^{(n)},\dots,R_1^{(n)})-\varphi_k(F(Z_{k+1}^{(n)}),\dots,F(Z_1^{(n)}))\}^2|\mathbf{Z}_{(\cdot)}^{(n)}]\\&\quad=o_{\mathrm{P}}(1)\end{aligned}\tag{4.1}$$

under $\mathcal{H}_{0;f}^{(n)}$, as $n\to\infty$. Once more, (4.1) automatically holds if, under $\mathcal{H}_{0;f}^{(n)}$,

$$\begin{aligned}&\mathrm{E}[\{a_k^{(n)}(\mathbf{N}^{(n)};R_{k+1}^{(n)},\dots,R_1^{(n)})-\varphi_k(F(Z_{k+1}^{(n)}),\dots,F(Z_1^{(n)}))\}^2|\mathbf{N}^{(n)}]\\&\quad=o_{\mathrm{P}}(1)\end{aligned}\tag{4.2}$$

as $n\to\infty$. We then have the following conditional asymptotic representation and asymptotic normality results, which are the serial counterpart of Lemma 3.1.

LEMMA 4.1. *Let $\varphi_k:(0,1)^{k+1}\to\mathbb{R}$ be a score-generating function for $a_k^{(n)}$. Then:*

(i) (Asymptotic representation) *under $\mathcal{H}_{0;f}^{(n)}$, as $n\to\infty$,*

$$S_k^{(n)}-\mathrm{E}[S_k^{(n)}|\mathbf{N}^{(n)}]=T_{\varphi_k;f;k}^{(n)}-\mathrm{E}[T_{\varphi_k;f;k}^{(n)}|\mathbf{Z}_{(\cdot)}^{(n)}]+o_{\mathrm{P}}(1/\sqrt{n}\,),\tag{4.3}$$

*where*

$$T_{\varphi_k;f;k}^{(n)}:=\frac{1}{n-k}\sum_{t=k+1}^{n}\varphi_k(F(Z_t^{(n)}),\dots,F(Z_{t-k}^{(n)}))$$

*and*

$$\begin{aligned}&\mathrm{E}[T_{\varphi_k;f;k}^{(n)}|\mathbf{Z}_{(\cdot)}^{(n)}]\\&\quad=[n(n-1)\cdots(n-k)]^{-1}\sum_{1\le t_1\ne\cdots\ne t_{k+1}\le n}\!\!\cdots\!\!\sum\varphi_k(F(Z_{t_1}^{(n)}),\dots,F(Z_{t_{k+1}}^{(n)}));\end{aligned}$$

(ii) (Asymptotic normality) *if, moreover, $0<\int_{(0,1)^{k+1}}|\varphi_k(u_{k+1},\dots,u_1)|^{2+\delta}du_1\cdots du_{k+1}<\infty$ for some $\delta>0$, then, under $\mathcal{H}_0^{(n)}$, as $n\to\infty$,*

$$\sqrt{n-k}(S_k^{(n)}-\mathrm{E}[S_k^{(n)}|\mathbf{N}^{(n)}])\xrightarrow{\mathcal{L}}\mathcal{N}(0,V^2),$$

*where, denoting by* $U_1, U_2, \ldots$ *an i.i.d. sequence of standard uniformly distributed random variables,*

$$V^2 := \mathrm{E}[\{\varphi_k^*(U_{k+1}, \ldots, U_1)\}^2] \\ + 2\sum_{j=1}^{k} \mathrm{E}[\varphi_k^*(U_{k+1}, \ldots, U_1)\varphi_k^*(U_{k+1+j}, \ldots, U_{1+j})] \tag{4.4}$$

*with, for* $u_1, \ldots, u_{k+1} \in (0,1)$,

$$\varphi_k^*(u_{k+1}, \ldots, u_1) \\ := \varphi_k(u_{k+1}, \ldots, u_1) \\ - \sum_{l=1}^{k+1} \mathrm{E}[\varphi_k(U_{k+1}, \ldots, U_1)|U_l = u_1] + k\mathrm{E}[\varphi_k(U_{k+1}, \ldots, U_1)].$$

PROOF. To prove part (i) of the lemma, we only need to show that, under $\mathcal{H}_{0;f}^{(n)}$, as $n \to \infty$, $\mathrm{E}[\{D_k^{(n)}\}^2|\mathbf{Z}_{(\cdot)}^{(n)}] = o_{\mathrm{P}}(1)$, where

$$D_k^{(n)} := \sqrt{n-k}\{(S_k^{(n)} - \mathrm{E}[S_k^{(n)}|\mathbf{N}^{(n)}]) - (T_{\varphi_k;f;k}^{(n)} - \mathrm{E}[T_{\varphi_k;f;k}^{(n)}|\mathbf{Z}_{(\cdot)}^{(n)}])\}.$$

Since the maximal invariant $(\mathbf{N}^{(n)}, \mathbf{R}^{(n)})$ depends on $\mathbf{Z}_{(\cdot)}^{(n)}$ only through $\mathbf{N}^{(n)}$, we actually have $\mathrm{E}[\{D_k^{(n)}\}^2|\mathbf{Z}_{(\cdot)}^{(n)}] = (n-k)\,\mathrm{Var}[S_k^{(n)} - T_{\varphi_k;f;k}^{(n)}|\mathbf{Z}_{(\cdot)}^{(n)}]$. Conditionally on $\mathbf{Z}_{(\cdot)}^{(n)}$ (and hence on $\mathbf{N}^{(n)}$), $S_k^{(n)} - T_{\varphi_k;f;k}^{(n)}$ is a linear serial *rank* statistic in the sense of Hallin, Ingenbleek and Puri [5]. Corollary 2 of Lemma 2, and Lemma 4 (Appendix 3) of that paper imply that there exists a constant $K$ (not depending on $n$) such that

$$\mathrm{E}[\{D_k^{(n)}\}^2|\mathbf{Z}_{(\cdot)}^{(n)}] \le \left(2k + 1 + \frac{K}{n-k}\right) \\ \times \mathrm{E}[\{a_k^{(n)}(\mathbf{N}^{(n)}; R_{k+1}^{(n)}, \ldots, R_1^{(n)}) \\ - \varphi_k(F(Z_{k+1}^{(n)}), \ldots, F(Z_1^{(n)}))\}^2|\mathbf{Z}_{(\cdot)}^{(n)}].$$

By (4.1), the last term converges to zero in probability under $\mathcal{H}_{0;f}^{(n)}$ as $n \to \infty$, which completes the proof of (4.3).

The asymptotic normality of $\sqrt{n-k}(T_{\varphi_k;f;k}^{(n)} - \mathrm{E}[T_{\varphi_k;f;k}^{(n)}|\mathbf{Z}_{(\cdot)}^{(n)}])$ [part (ii) of Lemma 4.1], hence also that of $\sqrt{n-k}(S_k^{(n)} - \mathrm{E}[S_k^{(n)}|\mathbf{N}^{(n)}])$, is also established in [5]. The special form of $V^2$ follows from Yoshihara's [20] central limit theorem for $U$-statistics under absolutely regular processes, which requires the $(2+\delta)$-integrability of the score-generating function $\varphi_k$. $\square$

Note that the right-hand side in (4.3) is exactly the same as in the asymptotic representation of the purely rank-based serial statistic

$$(n-k)^{-1}\sum_{t=k+1}^{n}\varphi_k\biggl(\frac{R_t^{(n)}}{n+1},\dots,\frac{R_{t-k}^{(n)}}{n+1}\biggr) - \mathrm{E}\biggl[(n-k)^{-1}\sum_{t=k+1}^{n}\varphi_k\biggl(\frac{R_t^{(n)}}{n+1},\dots,\frac{R_{t-k}^{(n)}}{n+1}\biggr)\biggr].$$

This remark, which is analogous to the remark made in the nonserial case just before the proof of Lemma 3.1, will play a crucial role in the proof of the asymptotic normality part of Proposition 4.1(ii).

4.2. *Exact and approximate scores.* As in the nonserial case, two types of scores—the exact and the approximate ones—are naturally associated with a given score-generating function. Define (referring to Section 3.2 for notation)

$$S_{\varphi_k;\mathrm{ex/appr}}^{(n)} := \frac{1}{n-k}\sum_{t=k+1}^{n} a_{\varphi_k;\mathrm{ex/appr}}^{(n)}(\mathbf{N}^{(n)};R_t^{(n)},R_{t-1}^{(n)},\dots,R_{t-k}^{(n)}),$$

where, for $(\eta,\nu)\in\{0,1,\dots,n\}^2$, $\nu\le n-\eta$ and $1\le i_1\ne i_2\ne\cdots\ne i_{k+1}\le n$,

$$\begin{aligned}a_{\varphi_k;\mathrm{ex}}^{(n)}&((\eta,\nu);i_1,\dots,i_{k+1})\\ &:=\mathrm{E}[\varphi_k(U_1^{(n)},\dots,U_{k+1}^{(n)})|N_{\mathbf{U};-}^{(n)}=\eta,N_{\mathbf{U};+}^{(n)}=\nu,\\ &\qquad\qquad R_{U_1}^{(n)}=i_1,\dots,R_{U_{k+1}}^{(n)}=i_{k+1}]\end{aligned}$$

and

$$\begin{aligned}a_{\varphi_k;\mathrm{appr}}^{(n)}&((\eta,\nu);i_1,\dots,i_{k+1})\\ &:=\varphi_k(\mathrm{E}[U_1^{(n)}|N_{\mathbf{U};-}^{(n)}=\eta,N_{\mathbf{U};+}^{(n)}=\nu,R_{U_1}^{(n)}=i_1],\\ &\qquad\dots,\mathrm{E}[U_{k+1}^{(n)}|N_{\mathbf{U};-}^{(n)}=\eta,N_{\mathbf{U};+}^{(n)}=\nu,R_{U_{k+1}}^{(n)}=i_{k+1}])\\ &=\varphi_k\biggl(I[i_1\le\eta]\biggl(\frac{i_1}{2(\eta+1)}\biggr)+I[i_1>n-\nu]\biggl(\frac{1}{2}+\frac{i_1-(n-\nu)}{2(\nu+1)}\biggr),\\ &\qquad\dots,I[i_{k+1}\le\eta]\biggl(\frac{i_{k+1}}{2(\eta+1)}\biggr)\\ &\qquad\qquad+I[i_{k+1}>n-\nu]\biggl(\frac{1}{2}+\frac{i_{k+1}-(n-\nu)}{2(\nu+1)}\biggr)\biggr).\end{aligned}$$

The following lemma provides sufficient conditions for $\varphi_k$ to be a score-generating function for $a_{\varphi_k;\mathrm{ex}}^{(n)}$ and for $a_{\varphi_k;\mathrm{appr}}^{(n)}$.

LEMMA 4.2. *Let $\varphi_k : (0,1)^{k+1} \longrightarrow \mathbb{R}$ be nonconstant and square-integrable. Then $\varphi_k$ is a score-generating function for $a^{(n)}_{\varphi_k;\mathrm{ex}}$. If, moreover, $\varphi_k$ is a linear combination of a finite number of square-integrable functions that are monotone in all their arguments, then $\varphi_k$ is also a score-generating function for $a^{(n)}_{\varphi_k;\mathrm{appr}}$.*

PROOF. The proof easily follows along the same lines as in the nonserial case and is left to the reader. □

4.3. *Unconditional asymptotic representation.* Lemma 4.1 was only an intermediate, conditional result; the following proposition provides the corresponding unconditional asymptotic representation and asymptotic normality. Define

$$\mu_{\varphi_k} := \int_{[0,1]^{k+1}} \varphi_k(u_0,\ldots,u_k)\,du_0\cdots du_k,$$

$$\mu^{(0)}_{\varphi_k} := \int_{[0,1/2]^{k+1}} \varphi_k(u_0,\ldots,u_k)\,du_0\cdots du_k,$$

$$\mu^{(k+1)}_{\varphi_k} := \int_{[1/2,1]^{k+1}} \varphi_k(u_0,\ldots,u_k)\,du_0\cdots du_k$$

and, for $\nu = 1,2,\ldots,k$,

$$\mu^{(\nu)}_{\varphi_k} := \sum_{0\le i_1<\cdots<i_\nu\le k}\cdots\sum \int_{(u_{i_1},\ldots,u_{i_\nu})\in[1/2,1]^\nu} \int_{(u_j,0\le j\le k,j\ne i_1,\ldots,i_\nu)\in[0,1/2]^{k+1-\nu}} \varphi_k(u_0,\ldots,u_k)\,du_0\cdots du_k.$$

PROPOSITION 4.1. *Let $\varphi_k$ be a nonconstant square-integrable score-generating function for $S^{(n)}_{\varphi_k;\mathrm{ex/appr}}$. Whenever approximate scores are considered, assume that $\varphi_k$ is a linear combination of square-integrable functions that are monotone in all their arguments. Then, writing $S^{(n)}_k$ for either $S^{(n)}_{\varphi_k;\mathrm{ex}}$ or $S^{(n)}_{\varphi_k;\mathrm{appr}}$:*

(i) (Asymptotic representation) *under $\mathcal{H}^{(n)}_{0;f}$, as $n\to\infty$,*

$$\begin{aligned} S^{(n)}_k - \mathrm{E}[S^{(n)}_k] \\ = T^{(n)}_{\varphi_k;f;k} - \mathrm{E}[T^{(n)}_{\varphi_k;f;k} | \mathbf{Z}^{(n)}_{(\cdot)}] \\ + 2^{k+1}[n(n-1)\cdots(n-k)]^{-1} \\ \times \Big\{ I[N^{(n)}_- \ge k+1] N^{(n)}_-(N^{(n)}_- - 1)\cdots(N^{(n)}_- - k)\mu^{(0)}_{\varphi_k} \end{aligned} \tag{4.5}$$

$$
\begin{aligned}
&+\sum_{\nu=1}^{k} I[k+1-\nu \le N_-^{(n)} \le n-\nu] \\
&\qquad \times N_-^{(n)}(N_-^{(n)}-1)\cdots(N_-^{(n)}-k+\nu) \\
&\qquad \times N_+^{(n)}(N_+^{(n)}-1)\cdots(N_+^{(n)}-\nu+1)\mu_{\varphi_k}^{(\nu)} \\
&+ I[N_+^{(n)} \ge k+1] N_+^{(n)}(N_+^{(n)}-1)\cdots(N_+^{(n)}-k)\mu_{\varphi_k}^{(k+1)} \Bigg\} \\
&- \mu_{\varphi_k} + o_{\mathrm{P}}(1/\sqrt{n});
\end{aligned}
$$

(ii) (Asymptotic normality) *if, moreover,* $\varphi_k$ *is* $(2+\delta)$*-integrable for some* $\delta > 0$*, then, under* $\mathcal{H}_0^{(n)}$*, as* $n \to \infty$*,*

$$
\begin{aligned}
&\sqrt{n-k}(S_k^{(n)} - \mathrm{E}[S_k^{(n)}]) \Big/ \sqrt{V^2 + (k+1)^2 \left[\mu_{\varphi_k} - 2\sum_{\nu=1}^{k+1} \nu \mu_{\varphi_k}^{(\nu)}/(k+1)\right]^2} \\
&\quad \xrightarrow{\mathcal{L}} \mathcal{N}(0,1),
\end{aligned} \tag{4.6}
$$

*with* $V^2$ *given in* (4.4).

When the score $\varphi_k$ is skew-symmetric with respect to $1/2$ [i.e., $\varphi_k(u_0, \ldots, u_i, \ldots, u_k) = -\varphi_k(u_0, \ldots, 1-u_i, \ldots, u_k)$ for all $i = 1, \ldots, k$], then $\mu_{\varphi_k} = \sum_{\nu=0}^{k+1} \mu_{\varphi_k}^{(\nu)} = 0$ with

$$
\mu_{\varphi_k}^{(\nu)} = \binom{k+1}{\nu} \int_{[1/2,1]^\nu \times [0,1/2]^{k+1-\nu}} \varphi_k(u_0, \ldots, u_k)\, du_0 \cdots du_k,
$$

so that

$$
\mu_{\varphi_k}^{(0)} = -\binom{k+1}{1}^{-1} \mu_{\varphi_k}^{(1)} = \binom{k+1}{2}^{-1} \mu_{\varphi_k}^{(2)} = -\binom{k+1}{3}^{-1} \mu_{\varphi_k}^{(3)} = \cdots.
$$

This and the fact that $N_-^{(n)}/n - \frac{1}{2} = O_{\mathrm{P}}(n^{-1/2})$ implies that the right-hand side of (4.5) reduces to $T_{\varphi_k;f;k}^{(n)} - \mathrm{E}[T_{\varphi_k;f;k}^{(n)} | \mathbf{Z}_{(\cdot)}^{(n)}] + o_{\mathrm{P}}(1)$. Hence, the conditional (4.3) and unconditional (4.5) representations of $S_k^{(n)} - \mathrm{E}[S_k^{(n)}]$ coincide.

Proof of Proposition 4.1. As in the nonserial case, we first prove the asymptotic representation result for exact scores. From the definition of exact scores, we obtain, for $S_k^{(n)} = S_{\varphi_k;\mathrm{ex}}^{(n)}$, writing $T_k^{(n)}$ for $T_{\varphi_k;f;k}^{(n)} := \frac{1}{n-k}\sum_{t=k+1}^{n} \varphi_k(F(Z_t^{(n)}), \ldots, F(Z_{t-k}^{(n)}))$,

$$
\mathrm{E}[S_k^{(n)} | \mathbf{N}^{(n)}]
$$

$$= \mathrm{E}[\mathrm{E}[T_k^{(n)}|\mathbf{R}^{(n)}, \mathbf{N}^{(n)}]|\mathbf{N}^{(n)}] = \mathrm{E}[T_k^{(n)}|\mathbf{N}^{(n)}]$$

$$= \mathrm{E}\left[\frac{1}{n-k}\sum_{t=k+1}^{n} \mathrm{E}[\varphi_k(F(Z_t^{(n)}), \ldots, F(Z_{t-k}^{(n)}))|\mathbf{N}^{(n)}, s_t^{(n)}, \ldots, s_{t-k}^{(n)}]|\mathbf{N}^{(n)}\right],$$

where

$$\begin{aligned} &\mathrm{E}[\varphi_k(F(Z_t^{(n)}), \ldots, F(Z_{t-k}^{(n)}))|\mathbf{N}^{(n)}, s_t^{(n)}, \ldots, s_{t-k}^{(n)}] \\ &\qquad = 2^{k+1}\int_{[0,1]^{k+1}} \varphi_k(u_0, \ldots, u_k) \\ &\qquad\qquad \times I[\operatorname{sign}(u_0 - \tfrac{1}{2}) = s_t^{(n)}, \\ &\qquad\qquad\qquad \ldots, \operatorname{sign}(u_k - \tfrac{1}{2}) = s_{t-k}^{(n)}]\, du_0 \cdots du_k. \end{aligned} \tag{4.7}$$

The asymptotic representation (4.5) (for exact scores) follows by combining (4.7) and part (i) of Lemma 4.1. Turning to approximate scores, it is sufficient for (4.5) to hold that

$$E^{(n)} := \{\mathrm{E}[S^{(n)}_{\varphi_k;\mathrm{appr}}|\mathbf{N}^{(n)}] - \mathrm{E}[S^{(n)}_{\varphi_k;\mathrm{appr}}]\} - \{\mathrm{E}[S^{(n)}_{\varphi_k;\mathrm{ex}}|\mathbf{N}^{(n)}] - \mathrm{E}[S^{(n)}_{\varphi_k;\mathrm{ex}}]\} \tag{4.8}$$

be $o_{\mathrm{P}}(1/\sqrt{n}\,)$. Note that

$$\begin{aligned} &\mathrm{E}[S^{(n)}_{\varphi_k;\mathrm{appr}}|\mathbf{N}^{(n)}] \\ &\quad = [n(n-1)\cdots(n-k)]^{-1} \\ &\qquad \times \sum_{1\le i_1\ne\cdots\ne i_{k+1}\le n}\cdots\sum \\ &\qquad\quad \varphi_k\biggl(I[i_1 \le N_-^{(n)}]\biggl(\frac{i_1}{2(N_-^{(n)}+1)}\biggr) + I[i_1 > N_-^{(n)}]\biggl(\frac{1}{2} + \frac{i_1 - N_-^{(n)}}{2(N_+^{(n)}+1)}\biggr), \\ &\qquad\qquad \ldots, I[i_{k+1} \le N_-^{(n)}]\biggl(\frac{i_{k+1}}{2(N_-^{(n)}+1)}\biggr) \\ &\qquad\qquad\qquad + I[i_{k+1} > N_-^{(n)}]\biggl(\frac{1}{2} + \frac{i_{k+1} - N_-^{(n)}}{2(N_+^{(n)}+1)}\biggr)\biggr). \end{aligned}$$

For notational simplicity, let us consider the case $k = 1$; the general case follows along the same ideas. For $k = 1$ we have

$$\begin{aligned} &\mathrm{E}[S^{(n)}_{\varphi_1;\mathrm{appr}}|\mathbf{N}^{(n)}] - \mathrm{E}[S^{(n)}_{\varphi_1;\mathrm{ex}}|\mathbf{N}^{(n)}] \\ &\quad = \frac{1}{n(n-1)}\Biggl\{\sum_{1\le i\ne j\le N_-^{(n)}}\sum \varphi_1\biggl(\frac{i}{2(N_-^{(n)}+1)}, \frac{j}{2(N_-^{(n)}+1)}\biggr) \end{aligned}$$

$$
\begin{aligned}
&+\sum_{i=1}^{N_-^{(n)}}\sum_{j=N_-^{(n)}+1}^{n}\varphi_1\left(\frac{i}{2(N_-^{(n)}+1)},\frac{1}{2}+\frac{j-N_-^{(n)}}{2(N_+^{(n)}+1)}\right)\\
&+\sum_{i=N_-^{(n)}+1}^{n}\sum_{j=1}^{N_-^{(n)}}\varphi_1\left(\frac{1}{2}+\frac{i-N_-^{(n)}}{2(N_+^{(n)}+1)},\frac{j}{2(N_-^{(n)}+1)}\right)\\
&+\sum\sum_{N_-^{(n)}+1\le i\ne j\le n}\varphi_1\left(\frac{1}{2}+\frac{i-N_-^{(n)}}{2(N_+^{(n)}+1)},\frac{1}{2}+\frac{j-N_-^{(n)}}{2(N_+^{(n)}+1)}\right)\Bigg\}\\
&-\frac{4}{n(n-1)}\{I[N_-^{(n)}\ge 2]N_-^{(n)}(N_-^{(n)}-1)\mu_{\varphi_1}^{(0)}\\
&\qquad+I[1\le N_-^{(n)}\le n-1]N_-^{(n)}N_+^{(n)}\mu_{\varphi_1}^{(1)}\\
&\qquad+I[N_+^{(n)}\ge 2]N_+^{(n)}(N_+^{(n)}-1)\mu_{\varphi_1}^{(2)}\}
\end{aligned}\tag{4.9}
$$

$$
\begin{aligned}
&=\frac{4N_-^{(n)}(N_-^{(n)}-1)^+}{n(n-1)}\Bigg\{\frac{(N_-^{(n)})^2}{N_-^{(n)}(N_-^{(n)}-1)}D^{--}_{N_-^{(n)},N_-^{(n)}}-\mu_{\varphi_1}^{--}\\
&\qquad-\frac{(N_-^{(n)})^2}{N_-^{(n)}(N_-^{(n)}-1)}\frac{1}{4(N_-^{(n)})^2}\\
&\qquad\times\sum_{i=1}^{N_-^{(n)}}\varphi_1\left(\frac{i}{2(N_-^{(n)}+1)},\frac{i}{2(N_-^{(n)}+1)}\right)\Bigg\}\\
&+\frac{4N_-^{(n)}N_+^{(n)}}{n(n-1)}\{D^{-+}_{N_-^{(n)},N_+^{(n)}}-\mu_{\varphi_1}^{-+}\}+\frac{4N_+^{(n)}N_-^{(n)}}{n(n-1)}\{D^{+-}_{N_+^{(n)},N_-^{(n)}}-\mu_{\varphi_1}^{+-}\}\\
&+\frac{4N_+^{(n)}(N_+^{(n)}-1)^+}{n(n-1)}\Bigg\{\frac{(N_+^{(n)})^2}{N_+^{(n)}(N_+^{(n)}-1)}D^{++}_{N_+^{(n)},N_+^{(n)}}-\mu_{\varphi_1}^{++}\\
&\qquad-\frac{(N_+^{(n)})^2}{N_+^{(n)}(N_+^{(n)}-1)}\frac{1}{4(N_+^{(n)})^2}\\
&\qquad\times\sum_{i=1}^{N_+^{(n)}}\varphi_1\left(\frac{1}{2}+\frac{i}{2(N_+^{(n)}+1)},\frac{1}{2}+\frac{i}{2(N_+^{(n)}+1)}\right)\Bigg\},
\end{aligned}
$$

where $x^+ := \max(0,x)$,

$$D^{--}_{\ell,m} := \frac{1}{4\ell m}\sum_{i=1}^{\ell}\sum_{j=1}^{m}\varphi_1\biggl(\frac{i}{2(\ell+1)},\frac{j}{2(m+1)}\biggr),$$

$$D^{-+}_{\ell,m} := \frac{1}{4\ell m}\sum_{i=1}^{\ell}\sum_{j=1}^{m}\varphi_1\biggl(\frac{i}{2(\ell+1)},\frac{1}{2}+\frac{j}{2(m+1)}\biggr),$$

$$D^{+-}_{\ell,m} := \frac{1}{4\ell m}\sum_{i=1}^{\ell}\sum_{j=1}^{m}\varphi_1\biggl(\frac{1}{2}+\frac{i}{2(\ell+1)},\frac{j}{2(m+1)}\biggr),$$

$$D^{++}_{\ell,m} := \frac{1}{4\ell m}\sum_{i=1}^{\ell}\sum_{j=1}^{m}\varphi_1\biggl(\frac{1}{2}+\frac{i}{2(\ell+1)},\frac{1}{2}+\frac{j}{2(m+1)}\biggr)$$

are Riemann sums for the integrals

$$\mu^{--}_{\varphi_1} := \int_0^{1/2}\int_0^{1/2}\varphi_1(u_0,u_1)\,du_0\,du_1,$$

$$\mu^{-+}_{\varphi_1} := \int_0^{1/2}\int_{1/2}^{1}\varphi_1(u_0,u_1)\,du_0\,du_1,$$

$$\mu^{+-}_{\varphi_1} := \int_{1/2}^{1}\int_0^{1/2}\varphi_1(u_0,u_1)\,du_0\,du_1,$$

$$\mu^{++}_{\varphi_1} := \int_{1/2}^{1}\int_{1/2}^{1}\varphi_1(u_0,u_1)\,du_0\,du_1,$$

respectively. Here again, due to the fact that $\varphi_1$ is square-integrable, the function $(u,v)\mapsto\varphi_1^*(u,v) := \varphi_1(u,v)I[u=v]$, $(u,v)\in[1/2,1]^2$, which vanishes except over the diagonal of the unit square, is integrable and has integral zero. Hence, $(1/4m^2)\sum_{i=1}^m\varphi_1^2(\frac{1}{2}+\frac{i}{2(m+1)},\frac{1}{2}+\frac{i}{2(m+1)})$, as a Riemann sum for the integral of $\varphi_1^*$ over $[1/2,1]^2$, is $o(1)$. Since

$$\biggl[\sum_{i=1}^{m}\varphi_1\biggl(\frac{1}{2}+\frac{i}{2(m+1)},\frac{1}{2}+\frac{i}{2(m+1)}\biggr)\biggr]^2 \leq m\sum_{i=1}^{m}\varphi_1^2\biggl(\frac{1}{2}+\frac{i}{2(m+1)},\frac{1}{2}+\frac{i}{2(m+1)}\biggr),$$

it follows that $(1/4m^2)\sum_{i=1}^m\varphi_1(\frac{1}{2}+\frac{i}{2(m+1)},\frac{1}{2}+\frac{i}{2(m+1)})$ is $o(1/\sqrt{m})$, as $m\to\infty$. A similar result holds for $(1/4m^2)\sum_{i=1}^m\varphi_1(\frac{i}{2(m+1)},\frac{i}{2(m+1)})$, as well as, of course, for any individual terms such as $(1/4m^2)\varphi_1(\frac{1}{2}+\frac{m}{2(m+1)},\frac{1}{2}+\frac{m}{2(m+1)})$.

Thus, (4.9) as $n \to \infty$ takes the form

$$\begin{aligned}
&\mathrm{E}[S^{(n)}_{\varphi_1;\mathrm{appr}} | \mathbf{N}^{(n)}] - \mathrm{E}[S^{(n)}_{\varphi_1;\mathrm{ex}} | \mathbf{N}^{(n)}] \\
&\quad = \frac{4N^{(n)}_{-}(N^{(n)}_{-} - 1)^{+}}{n(n-1)} [D^{--}_{N^{(n)}_{-}, N^{(n)}_{-}} - \mu^{--}_{\varphi_1}] \\
&\qquad + \frac{4N^{(n)}_{-} N^{(n)}_{+}}{n(n-1)} [D^{-+}_{N^{(n)}_{-}, N^{(n)}_{+}} - \mu^{-+}_{\varphi_1}] \\
&\qquad + \frac{4N^{(n)}_{+} N^{(n)}_{-}}{n(n-1)} [D^{+-}_{N^{(n)}_{+}, N^{(n)}_{-}} - \mu^{+-}_{\varphi_1}] \\
&\qquad + \frac{4N^{(n)}_{+}(N^{(n)}_{+} - 1)^{+}}{n(n-1)} [D^{++}_{N^{(n)}_{+}, N^{(n)}_{+}} - \mu^{++}_{\varphi_1}] + o_{\mathrm{P}}(1/\sqrt{n}).
\end{aligned}$$

Considering the difference $D^{++}_{m,m} - \mu^{++}_{\varphi_1}$, we have

$$\begin{aligned}
D^{++}_{m,m} - \mu^{++}_{\varphi_1} &= \frac{1}{4m^2} \sum_{i=1}^{m} \sum_{j=1}^{m} \varphi_1 \left( \frac{1}{2} + \frac{i}{2(m+1)}, \frac{1}{2} + \frac{j}{2(m+1)} \right) \\
&\quad - \iint_{[1/2,1]^2} \varphi_1(u_0, u_1)\, du_0\, du_1 \\
&= \frac{1}{4m^2} \sum_{i=1}^{m-1} \sum_{j=1}^{m-1} \varphi_1 \left( \frac{1}{2} + \frac{i}{2(m+1)}, \frac{1}{2} + \frac{j}{2(m+1)} \right) \\
&\quad - \iint_{[1/2,1]^2} \varphi_1(u_0, u_1)\, du_0\, du_1 + o(1/\sqrt{m}),
\end{aligned} \tag{4.10}$$

because, in view of the same argument as above, the two first sums in (4.10) are $o(1/\sqrt{m})$. As in the proof of Proposition 3.1, due to the fact that $\varphi_1$ can be assumed to be nondecreasing in its two arguments, the sum that appears in this latter expression is composed between the two Darboux sums

$$\underline{D}^{++}_{m,m} := \frac{1}{4m^2} \sum_{i=1}^{m-1} \sum_{j=1}^{m-1} \varphi_1 \left( \frac{1}{2} + \frac{i-1}{2m}, \frac{1}{2} + \frac{j-1}{2m} \right)$$

and

$$\bar{D}^{++}_{m,m} := \frac{1}{4m^2} \sum_{i=1}^{m-1} \sum_{j=1}^{m-1} \varphi_1 \left( \frac{1}{2} + \frac{i}{2m}, \frac{1}{2} + \frac{j}{2m} \right).$$

These Darboux sums also converge to the integral $\iint_{[1/2,1]^2} \varphi_1(u_0, u_1)\, du_0\, du_1$

and

$$\bar{D}^{++}_{m,m} - \underline{D}^{++}_{m,m} = \frac{1}{4m^2}\left[\varphi_1\left(\frac{1}{2}+\frac{m-1}{2m}, \frac{1}{2}+\frac{m-1}{2m}\right) - \varphi_1\left(\frac{1}{2},\frac{1}{2}\right)\right];$$

the same argument still implies that this difference, hence also $D^{++}_{m,m} - \mu^{++}_{\varphi_1}$, is $o(1/\sqrt{m})$. The other three quantities of the same type can be treated similarly. Uniform integrability and the fact that $N^{(n)}_{\pm}$ are $O_{\mathrm{P}}(n)$, as in the proof of Proposition 3.1, complete the proof that (4.8) is indeed $o_{\mathrm{P}}(1/\sqrt{n})$.

To conclude, we now prove the asymptotic normality result. Denote by $\Pi_{k+1}$ the set of permutations $\pi$ of $\{1,\ldots,k+1\}$. Then

$$\begin{aligned}
&\mathrm{E}[S^{(n)}_{\varphi_k;\mathrm{ex}}|\mathbf{N}^{(n)}] \\
&\quad = \binom{n}{k+1}^{-1} \sum\cdots\sum_{1\le t_1<\cdots<t_{k+1}\le n} \Biggl\{\sum_{\nu=0}^{k+1}\frac{1}{(k+1)!}2^{k+1}\mu^{(\nu)}_{\varphi_k} \\
&\qquad\qquad \times \sum_{\pi\in\Pi_{k+1}} I[s^{(n)}_{t_{\pi(1)}}=1,\ldots,s^{(n)}_{t_{\pi(\nu)}}=1, \\
&\qquad\qquad\qquad s^{(n)}_{t_{\pi(\nu+1)}}=-1,\ldots,s^{(n)}_{t_{\pi(k+1)}}=-1]\Biggr\};
\end{aligned}$$

hence $\mathrm{E}[S^{(n)}_{\varphi_k;\mathrm{ex}}|\mathbf{N}^{(n)}]$ is a $U$-statistic computed from the $n$-tuple $Z^{(n)}_1,\ldots,Z^{(n)}_n$ with kernel

$$\begin{aligned}
&h_k(z_1,\ldots,z_{k+1}) \\
&\quad = \sum_{\nu=0}^{k+1}\frac{2^{k+1}\mu^{(\nu)}_{\varphi_k}}{(k+1)!} \\
&\qquad \times \sum_{\pi\in\Pi_{k+1}} I[z_{\pi(1)}>0,\ldots,z_{\pi(\nu)}>0, z_{\pi(\nu+1)}\le 0,\ldots,z_{\pi(k+1)}\le 0].
\end{aligned}$$

Routine calculation yields, under $\mathcal{H}^{(n)}_{0;f}$,

$$\begin{aligned}
&\mathrm{E}[h_k(Z^{(n)}_1,\ldots,Z^{(n)}_{k+1})|Z^{(n)}_1] \\
&\quad = 2I[Z^{(n)}_1>0]\sum_{\nu=0}^{k+1}\frac{\nu}{k+1}\mu^{(\nu)}_{\varphi_k} + 2I[Z^{(n)}_1\le 0]\sum_{\nu=0}^{k+1}\frac{k+1-\nu}{k+1}\mu^{(\nu)}_{\varphi_k}
\end{aligned}$$

and

$$\mathrm{Var}(\mathrm{E}[h_k(Z^{(n)}_1,\ldots,Z^{(n)}_{k+1})|Z^{(n)}_1]) = \left\{\mu_{\varphi_k} - 2\sum_{\nu=1}^{k+1}\frac{\nu}{k+1}\mu^{(\nu)}_{\varphi_k}\right\}^2,$$

which is strictly positive. Classical results on $U$-statistics (see, e.g., [19]) then imply that, under $\mathcal{H}_{0;f}^{(n)}$, as $n \to \infty$,

$$(n-k)^{1/2}(\mathrm{E}[S_{\varphi_k;\mathrm{ex}}^{(n)}|\mathbf{N}^{(n)}] - \mathrm{E}[S_{\varphi_k;\mathrm{ex}}^{(n)}])$$
$$\xrightarrow{\mathcal{L}} \mathcal{N}\left(0, (k+1)^2\left\{\mu_{\varphi_k} - 2\sum_{\nu=1}^{k+1}\frac{\nu}{k+1}\mu_{\varphi_k}^{(\nu)}\right\}^2\right).$$

The same argument as in the nonserial case can be invoked to establish the asymptotic independence of the right-hand side in the conditional asymptotic representation (4.3) and $(n-k)^{1/2}(\mathrm{E}[S_{\varphi_k;\mathrm{ex}}^{(n)}|\mathbf{N}^{(n)}] - \mathrm{E}[S_{\varphi_k;\mathrm{ex}}^{(n)}])$. The result follows. □

4.4. *Example*: *first-order median moving average.* The central sequence (2.5) under $\mathrm{P}_{f;\theta}^{(n)}$ clearly [central sequences are always defined up to $o_\mathrm{P}(1)$ quantities] can be rewritten as

$$\Delta_f^{(n)} = \sqrt{n-1}\,r_{f;1}^{(n)} + o_\mathrm{P}(1),$$

where, defining $\varphi_f(u) := \frac{-f'}{f}(F^{-1}(u))$ and $\psi_f(u) := F^{-1}(u)$, $u \in (0,1)$,

$$r_{f;1}^{(n)} := \frac{1}{n-1}\sum_{t=2}^{n}\varphi_f(F(Z_t^{(n)}))\psi_f(F(Z_{t-1}^{(n)})),$$

a measure of serial dependence associated with $f$. With this notation, it clearly appears that $r_{f;1}^{(n)}$ is a particular case [letting $k = 1$ and $\varphi_1(u_0, u_1) := \varphi_f(u_0)\psi_f(u_1)$] of the statistic $T_{\varphi_k;f;k}^{(n)}$ considered in Lemma 4.1.

Define the serial linear sign-and-rank autocorrelation statistic of order 1 (based on exact scores) as $\underset{\sim}{r}_{f;1;\mathrm{ex}}^{(n)*} := \mathrm{E}[r_{f;1}^{(n)}|\mathbf{N}^{(n)}, \mathbf{R}^{(n)}]$. Proposition 4.1 implies that, under $\mathrm{P}_{f;\theta}^{(n)}$, as $n \to \infty$,

$$\underset{\sim}{r}_{f;1;\mathrm{ex}}^{(n)*} = r_{f;1}^{(n)*} + o_\mathrm{P}(1/\sqrt{n})$$

with

$$\begin{aligned} r_{f;1}^{(n)*} &= r_{f;1}^{(n)} - \mathrm{E}[r_{f;1}^{(n)}|\mathbf{Z}_{(\cdot)}^{(n)}] \\ &\quad + \frac{2^2}{n(n-1)}\biggl\{I[N_-^{(n)} \ge 2]N_-^{(n)}(N_-^{(n)}-1)(-f(0))\int_{-\infty}^0 xf(x)\,dx \\ &\qquad + I[1 \le N_-^{(n)} \le n-1]N_-^{(n)}N_+^{(n)} \\ &\qquad \times \left(-f(0)\int_0^\infty xf(x)\,dx + f(0)\int_{-\infty}^0 xf(x)\,dx\right) \end{aligned}$$

$$+ I[N_+^{(n)} \geq 2] f(0) \int_0^\infty x f(x)\, dx \Big\}$$

$$+ o_{\mathrm{P}}(1/\sqrt{n}\,)$$

$$= r_{f;1}^{(n)} - \mathrm{E}[r_{f;1}^{(n)} | \mathbf{Z}_{(\cdot)}^{(n)}] + 2f(0)\mu_f \frac{N_+^{(n)} - N_-^{(n)}}{n} + o_{\mathrm{P}}(1/\sqrt{n}\,)$$

$$= \frac{1}{n-1} \sum_{t=2}^{n} \frac{-f'}{f}(Z_t^{(n)}) Z_{t-1}^{(n)} - \frac{1}{n(n-1)} \sum\sum_{1 \leq t_1 \neq t_2 \leq n} \frac{-f'}{f}(Z_{t_1}^{(n)}) Z_{t_2}^{(n)}$$

$$+ 2f(0)\mu_f \frac{N_+^{(n)} - N_-^{(n)}}{n} + o_{\mathrm{P}}(1/\sqrt{n}\,)$$

$$= \frac{1}{n-1} \sum_{t=2}^{n} \frac{-f'}{f}(Z_t^{(n)})(Z_{t-1}^{(n)} - \mu_f) + 2f(0)\mu_f \frac{N_+^{(n)} - N_-^{(n)}}{n} + o_{\mathrm{P}}(1/\sqrt{n}\,).$$

Letting

$$\underset{\sim}{r}_{f;1;\mathrm{appr}}^{(n)} := \frac{1}{n-1} \sum_{t=2}^{n} \varphi_f(\underset{\sim}{R}_t^{(n)}) \psi_f(\underset{\sim}{R}_{t-1}^{(n)}) \tag{4.11}$$

with $\underset{\sim}{R}_t^{(n)}$ given in (3.22), the approximate score counterpart of $\underset{\sim}{r}_{f;1;\mathrm{ex}}^{(n)*}$ is, in view of Lemmas 4.1 and 4.2,

$$\underset{\sim}{r}_{f;1;\mathrm{appr}}^{(n)*} = \underset{\sim}{r}_{f;1;\mathrm{appr}}^{(n)} - \mathrm{E}[\underset{\sim}{r}_{f;1;\mathrm{appr}}^{(n)} | \mathbf{N}^{(n)}] + 2f(0)\mu_f \frac{N_+^{(n)} - N_-^{(n)}}{n}$$

$$= \frac{1}{n-1} \sum_{t=2}^{n} \varphi_f(\underset{\sim}{R}_t^{(n)}) \psi_f(\underset{\sim}{R}_{t-1}^{(n)}) \tag{4.12}$$

$$- \frac{1}{n(n-1)} \sum\sum_{1 \leq t_1 \neq t_2 \leq n} \varphi_f(\underset{\sim}{R}_{t_1}^{(n)}) \psi_f(\underset{\sim}{R}_{t_2}^{(n)}) + 2f(0)\mu_f \frac{N_+^{(n)} - N_-^{(n)}}{n}.$$

In conclusion, defining $\Delta_f^{(n)*} := \sqrt{n-1} r_{f;1}^{(n)*}$ and $\underset{\sim}{\Delta}_{f;\mathrm{ex/appr}}^{(n)*} := \sqrt{n-1} \times \underset{\sim}{r}_{f;1;\mathrm{ex/appr}}^{(n)*}$, we obtain

$$\Delta_f^{(n)*} = \underset{\sim}{\Delta}_{f;\mathrm{ex}}^{(n)*} + o_{\mathrm{P}}(1) = \underset{\sim}{\Delta}_{f;\mathrm{appr}}^{(n)*} + o_{\mathrm{P}}(1)$$

under $\mathcal{H}_{0;f}^{(n)}$, as $n \to \infty$. Using standard arguments, one easily verifies that $\Delta_f^{(n)*}$, $\underset{\sim}{\Delta}_{f;\mathrm{ex}}^{(n)*}$ and $\underset{\sim}{\Delta}_{f;\mathrm{appr}}^{(n)*}$ are indeed three versions of the semiparametrically efficient central sequence for $\theta$ in the model $\mathcal{E}_0^{(n)}$; again, the sign-and-rank-based $\underset{\sim}{\Delta}_{f;\mathrm{appr}}^{(n)*}$ can be used to perform semiparametrically efficient inference (tests, estimation, etc.) for the MA(1) coefficient $\theta$; see, for example, Section 11.9 of [16].

**5. Numerical study.** The finite-sample performance of the proposed test statistics has been studied in the context of the first-order moving average model of the example in Section 4.4. More precisely, we generated 1000 replications of each of the MA(1) processes characterized by equation (2.4) with parameter values $\theta = \pm 0.3$, $\pm 0.25$, $\pm 0.20$, $\pm 0.15$, $\pm 0.10$, $\pm 0.05$ and 0, and the following asymmetric innovation densities:

(a) $f(z) := f_{t_1} I[z \le 0] + f_{\mathcal{N}(0,1)}(z) I[z > 0]$, where $f_{t_1}$ stands for the Cauchy density and $f_{\mathcal{N}(0,1)}$ for the standard normal one;
(b) $f(z) := f_{t_5} I[z \le 0] + f_{\mathcal{N}(0,1)}(z) I[z > 0]$, where $f_{t_5}$ stands for the Student density with 5 degrees of freedom;
(c) $f := f_{\mathcal{N}_{\lambda=-10}}$ (the skew normal density with skewness $\lambda = -10$; see [1]), duly shifted and rescaled to have zero median and unit variance;
(d) $f := f_{\mathcal{N}_{\lambda=-20}}$ (the skew normal density with skewness $\lambda = -20$), duly shifted and rescaled to have zero median and unit variance;
(e) $f := 0.5 f_{\mathcal{N}(0,1)} + 0.5 f_{\mathcal{N}(-5,2)}$ (a mixed-normal density), duly shifted and rescaled to have zero median and unit variance;
(f) $f := 0.75 f_{\mathcal{N}(0,1)} + 0.25 f_{\mathcal{N}(-5,1)}$ (a mixed-normal density), duly shifted and rescaled to have zero median and unit variance.

For each replication, randomness (namely, $\theta = 0$) has been tested against first-order moving average dependence (two-sided test), based on the asymptotically normal distribution of:

(i) the ordinary first-order autocorrelation coefficient

$$r_1^{(n)} := (n-1)^{-1} \sum_{t=2}^{n} (Z_t - \bar{Z}^{(n)})(Z_{t-1} - \bar{Z}^{(n)}) \Big/ n^{-1} \sum_{t=1}^{n} (Z_t - \bar{Z}^{(n)})^2;$$

(ii) the "traditional" first-order van der Waerden rank autocorrelation coefficient

$$\underset{\sim}{r}_{\mathrm{vdW};1}^{(n)} := \Biggl\{ (n-1)^{-1} \sum_{t=2}^{n} \Phi^{-1}\biggl(\frac{R_t^{(n)}}{n+1}\biggr) \Phi^{-1}\biggl(\frac{R_{t-1}^{(n)}}{n+1}\biggr) - [n(n-1)]^{-1} \sum_{1 \le i \ne j \le n} \Phi^{-1}\biggl(\frac{i}{n+1}\biggr) \Phi^{-1}\biggl(\frac{j}{n+1}\biggr) \Biggr\} \Big/ \sigma_{\mathrm{vdW};1}^{(n)},$$

where $\Phi$ stands for the standard normal distribution function and $\sigma_{\mathrm{vdW};1}^{(n)}$ stands for the exact standardizing constant (see, e.g., [8]);

(iii) the "traditional" first-order Wilcoxon rank autocorrelation coefficient

$$\underset{\sim}{r}^{(n)}_{W;1} := \Biggl\{(n-1)^{-1}\sum_{t=2}^{n}\varphi_{\log}\biggl(\frac{R_t^{(n)}}{n+1}\biggr)\psi_{\log}\biggl(\frac{R_{t-1}^{(n)}}{n+1}\biggr) - [n(n-1)]^{-1}\sum\sum_{1\le i\ne j\le n}\varphi_{\log}\biggl(\frac{i}{n+1}\biggr)\psi_{\log}\biggl(\frac{j}{n+1}\biggr)\Biggr\}\Big/\sigma^{(n)}_{\mathrm{W};1},$$

with $\varphi_{\log}(u) := 2u-1$ and $\psi_{\log}(u) := \ln(\frac{u}{1-u})$, $u \in (0,1)$ ($\psi_{\log}$ is proportional to the inverse of the logistic distribution function); $\sigma^{(n)}_{\mathrm{W};1}$ stands for the exact standardizing constant (see, e.g., [8]);

(iv) the "traditional" first-order Laplace rank autocorrelation coefficient

$$\underset{\sim}{r}^{(n)}_{\mathrm{L};1} := \Biggl\{(n-1)^{-1}\sum_{t=2}^{n}\varphi_{\exp}\biggl(\frac{R_t^{(n)}}{n+1}\biggr)\psi_{\exp}\biggl(\frac{R_{t-1}^{(n)}}{n+1}\biggr) - [n(n-1)]^{-1}\sum\sum_{1\le i\ne j\le n}\varphi_{\exp}\biggl(\frac{i}{n+1}\biggr)\psi_{\exp}\biggl(\frac{j}{n+1}\biggr)\Biggr\}\Big/\sigma^{(n)}_{\mathrm{L};1},$$

with $\varphi_{\exp}(u) := \mathrm{sign}(2u-1)$ and

$$\psi_{\exp}(u) := \ln(2u)I[u \le 0.5] - \ln 2(1-u)I[u > 0.5], \qquad u \in (0,1)$$

($\psi_{\exp}$ is proportional to the inverse of the double-exponential distribution function); $\sigma^{(n)}_{\mathrm{L};1}$ stands for the exact standardizing constant (see, e.g., [8]);

(v) the first-order sign-and-rank autocorrelation coefficient $\underset{\sim}{r}^{(n)*}_{\mathrm{W/vdW};1}$ defined in (4.12), with the approximate scores $\varphi_f(u) = \frac{1}{\gamma}\varphi_{\log}(u)I[u \le 0.5] + \phi^{-1}(u)I[u > 0.5]$ and $\psi_f(u) = \gamma\psi_{\log}(u)I[u \le 0.5] + \phi^{-1}(u)I[u > 0.5]$ associated with a density $f(z) := \frac{1}{\gamma}\frac{\exp(z/\gamma)}{(1+\exp(z/\gamma))^2}I[z \le 0] + f_{\mathcal{N}(0,1)}(z)I[z > 0]$ (with $\gamma := \sqrt{\pi/8}$) that is logistic on the negative half-line and standard normal on the positive half-line (yielding Wilcoxon scores for the negative residuals and van der Waerden scores for the positive ones);

(vi) the first-order sign-and-rank autocorrelation coefficient $\underset{\sim}{r}^{(n)*}_{L/\mathrm{vdW};1}$ defined in (4.12), with the approximate scores

$$\varphi_f(u) = -\frac{1}{\gamma}I[u \le 0.5] + \phi^{-1}(u)I[u > 0.5]$$

and

$$\psi_f(u) = \gamma\psi_{\exp}(u)I[u \le 0.5] + \phi^{-1}(u)I[u > 0.5]$$

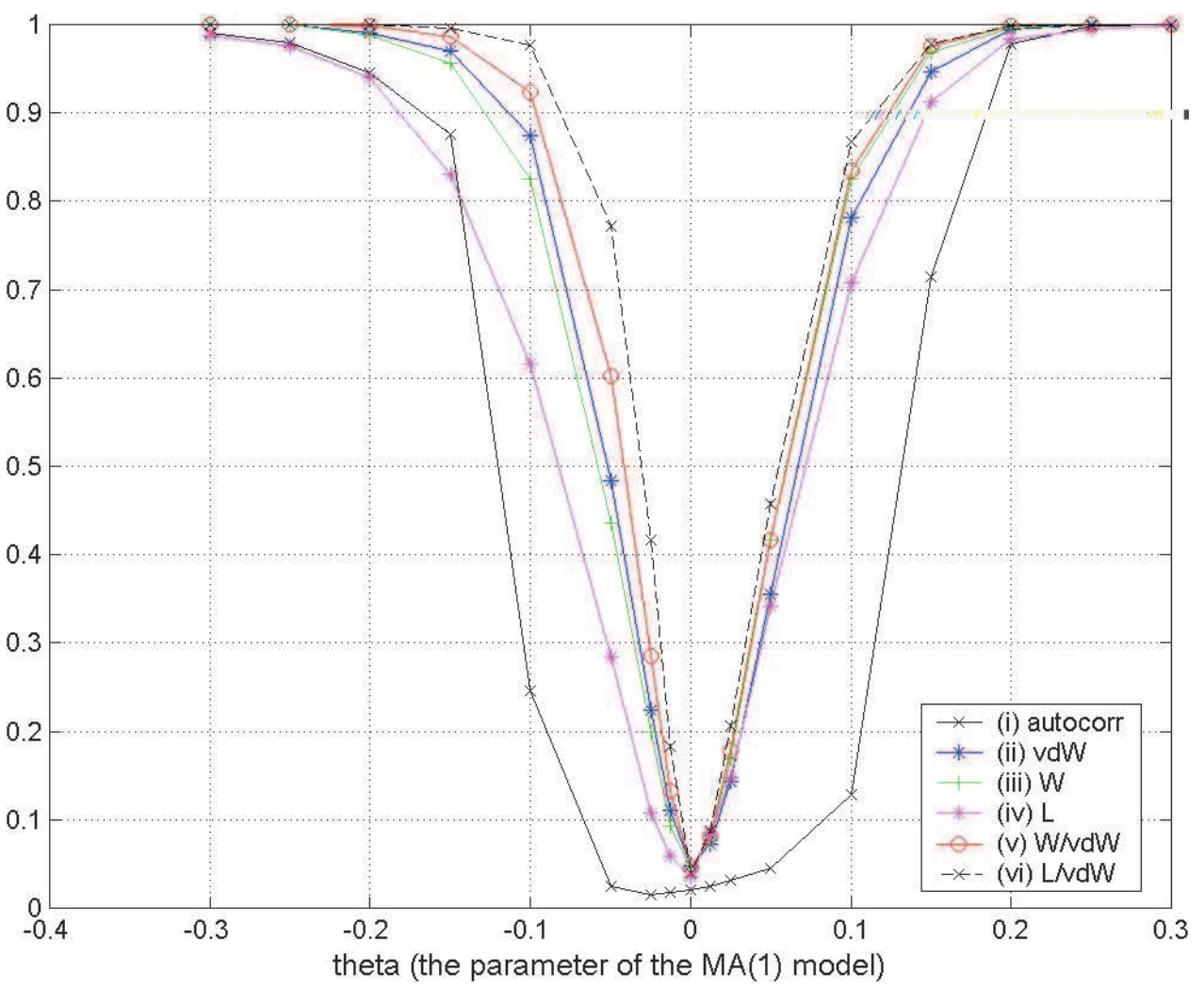

FIG. 1. *Empirical power, under Cauchy/standard normal innovations* (a), *of various parametric, rank and sign-and-rank tests for randomness against first-order MA dependence* [*based on the test statistics* (i)–(vi)]. *The series length is* $n = 250$; 1000 *replications were performed.*

associated with a density $f(z) := \frac{1}{2\gamma}\exp(z/\gamma)I[z \le 0] + f_{\mathcal{N}(0,1)}(z)I[z > 0]$ (with $\gamma = \sqrt{\pi/2}$) that is double-exponential on the negative half-line, and standard normal on the positive half-line (yielding Laplace scores for the negative residuals and van der Waerden scores for the positive ones).

The results of these simulations (series length $n = 250$; number of replications 1000) are summarized in Figures 1–6, where the graphs of the empirical power functions associated with testing procedures (i)–(vi) are plotted against $\theta$.

These graphs speak for themselves and need little comment. They all clearly demonstrate the superiority, under asymmetric densities, of the sign-and-rank methods over both their classical Gaussian and traditional rank-based competitors. The more skewed the underlying density, the more significant the improvement. For instance, in Figure 1 [Cauchy/Normal density (a)] the percentage of rejection at $\theta = -0.05$, which is only 0.0240 for the traditional correlogram-based tests (a severely biased test, thus), is as high as 0.7720 for the sign-and-rank Laplace/van der Waerden tests (vi). At $\theta = -0.10$, the corresponding figures are 0.2460 for the correlogram-based tests, but 0.9770 for the Laplace/van der Waerden ones. Of course, the performance of the parametric correlogram method in this case is particularly poor, due to the absence of finite moments, but the superiority of the sign-and-rank-based methods over their "purely rank-based" competitors remains quite substantial (at $\theta = -0.05$ and $\theta = -0.10$, Wilcoxon tests only yield empirical powers 0.4360 and 0.8250). Quite understandably, this superiority of the sign-and-rank methods over their competitors fades away under moderately skewed densities (see Figure 2, where it is less pronounced than in Figure 1), but it remains extremely substantial in Figures 4–6.

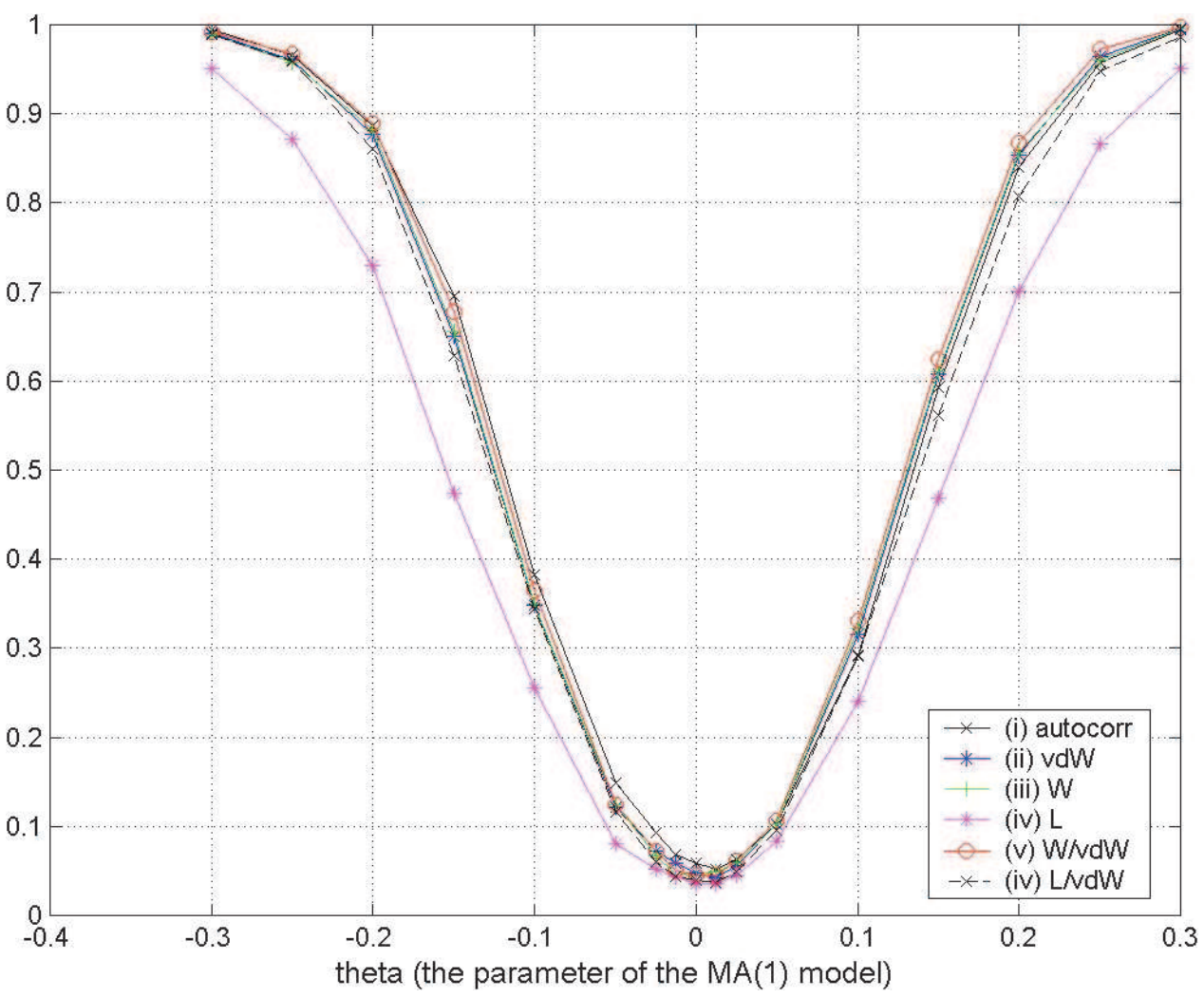

Fig. 2. *Empirical power, under Student (*5 *d.f.)/standard normal innovations* (b), *of various parametric, rank and sign-and-rank tests for randomness against first-order MA dependence* [*based on the test statistics* (i)–(vi)]. *The series length is* $n = 250$*;* 1000 *replications were performed.*

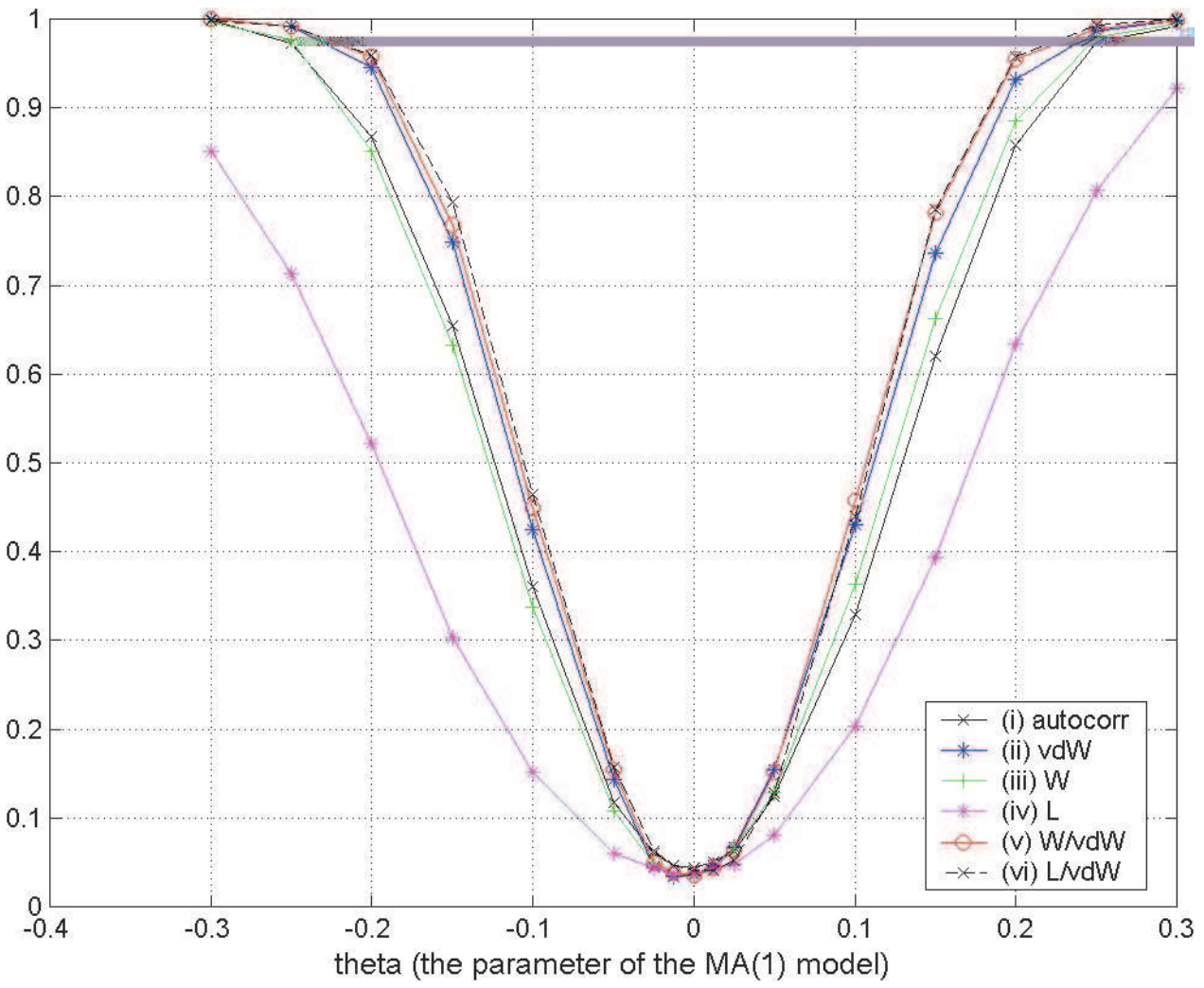

Fig. 3. *Empirical power, under skew-normal* ($\lambda = -10$) *innovations* (c), *of various parametric, rank and sign-and-rank tests for randomness against first-order MA dependence* [*based on the test statistics* (i)–(vi)]. *The series length is* $n = 250$*;* 1000 *replications were performed.*

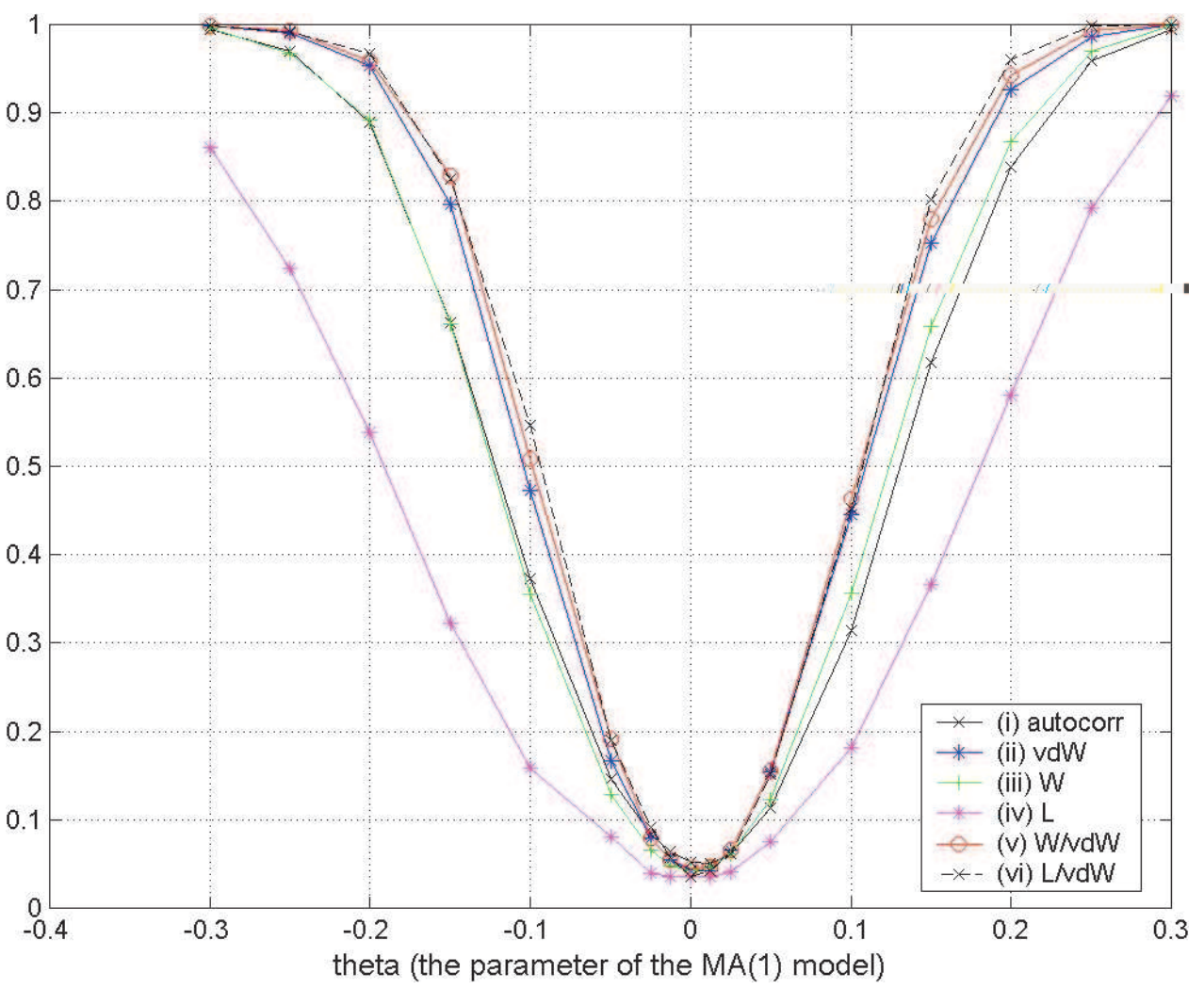

FIG. 4. *Empirical power, under under skew-normal ($\lambda = -20$) innovations* (d), *of various parametric, rank and sign-and-rank tests for randomness against first-order MA dependence* [*based on the test statistics* (i)–(vi)]. *The series length is* $n = 250$; 1000 *replications were performed.*

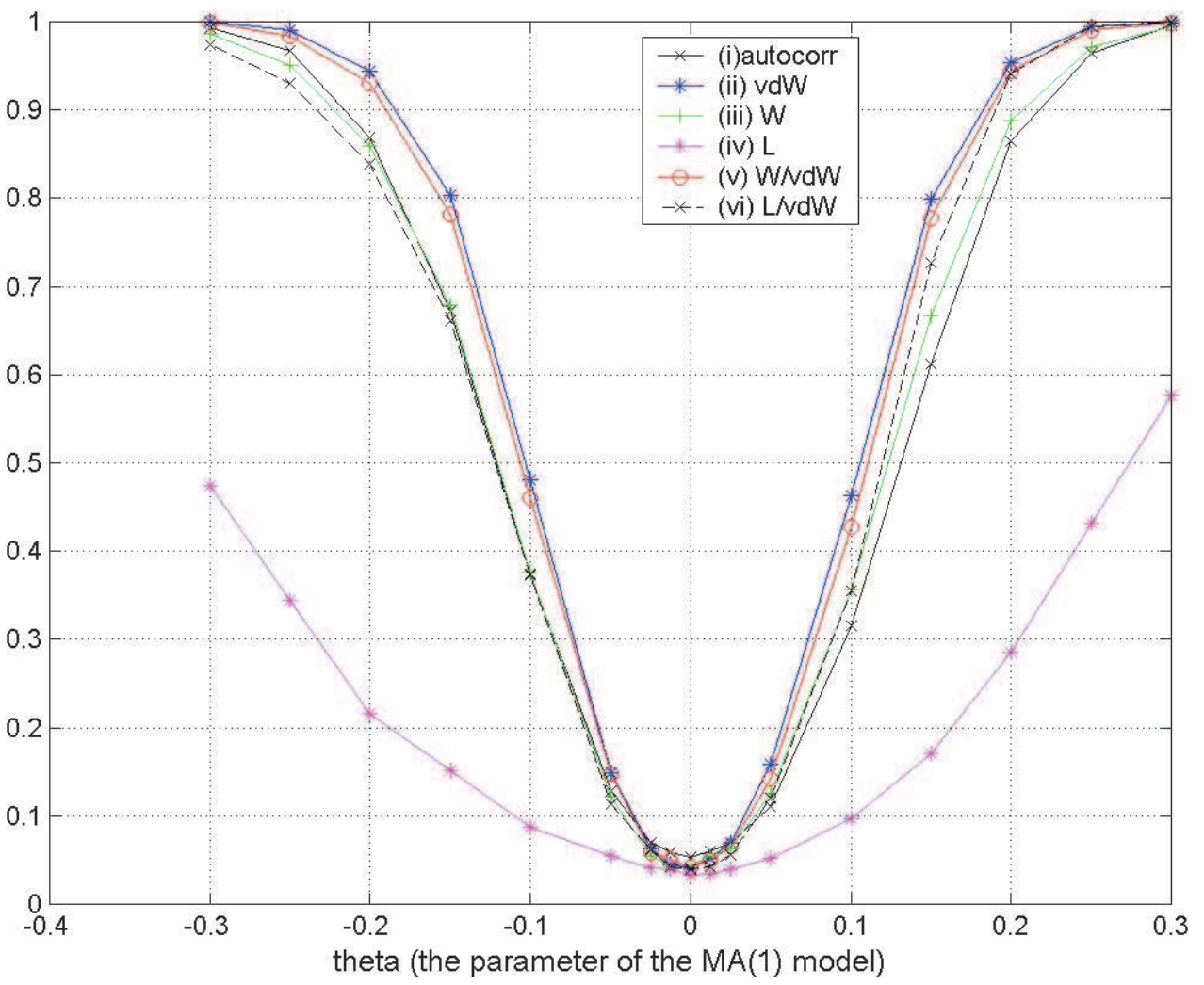

FIG. 5. *Empirical power, under mixed normal innovations* (e), *of various parametric, rank and sign-and-rank tests for randomness against first-order MA dependence* [*based on the test statistics* (i)–(vi)]. *The series length is* $n = 250$; 1000 *replications were performed.*

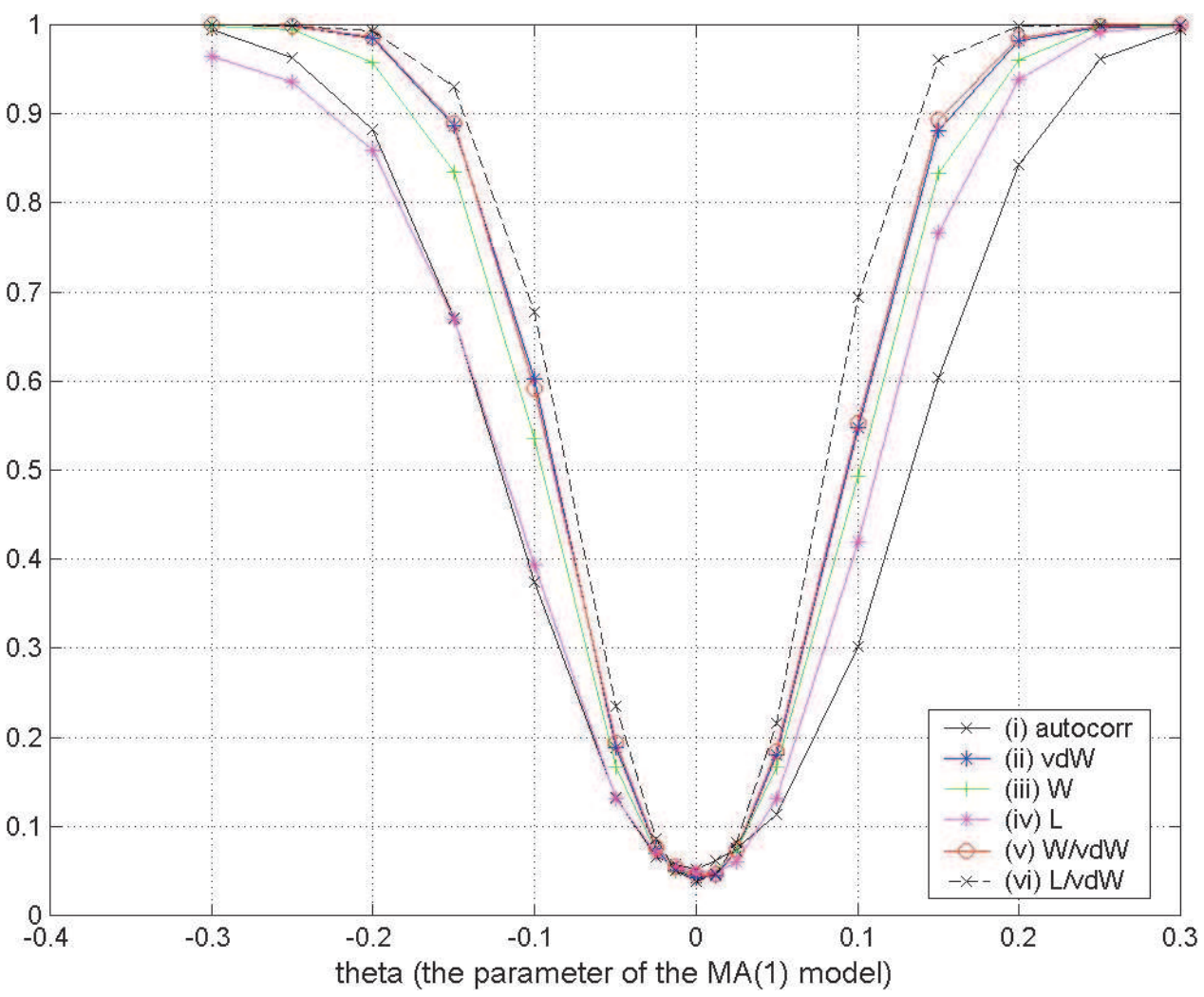

FIG. 6. *Empirical power, under mixed normal innovations* (f), *of various parametric, rank and sign-and-rank tests for randomness against first-order MA dependence* [*based on the test statistics* (i)–(vi)]. *The series length is* $n = 250$; 1000 *replications were performed.*

**Acknowledgments.** The authors would like to thank two referees and an Associate Editor for their pertinent remarks on the original manuscript and for their criticisms, which led to a complete rewriting of the paper.

M. Hallin
Institut de Statistique et
de Recherche Opérationnelle
ECARES
and
Département de Mathématique
Université Libre de Bruxelles
CP 210
B-1050 Brussels
Belgium
E-mail: mhallin@ulb.ac.be

C. Vermandele
Institut de Statistique et
de Recherche Opérationnelle
and Faculté des Sciences Sociales,
Politiques et Economiques
Université Libre de Bruxelles
CP 124
B-1050 Brussels
Belgium
E-mail: vermande@ulb.ac.be

B. Werker
Econometrics and Finance Group
CentER
Tilburg University
P.O. Box 90153
5000 LE Tilburg
The Netherlands
E-mail: B.J.M.Werker@uvt.nl